\documentclass[a4paper,fleqn]{cas-sc}
\usepackage{amsmath}
\usepackage{ifpdf}
\usepackage{indentfirst}
\usepackage{hyperref}
\usepackage{algorithm}
\usepackage{algorithmicx}  
\usepackage{algpseudocode}
\usepackage{subcaption}
\usepackage{booktabs}
\usepackage{multirow}
\usepackage{float}

\numberwithin{equation}{section}

\usepackage{graphicx,cite,amssymb,lineno,color,amsfonts}

\newtheorem{Theorem}{Theorem}[section]
\newtheorem{Lemma}{Lemma}[section]

\newtheorem{Definition}{Definition}[section]

\makeatletter
\def\elsartstyle{%
    \def\normalsize{\@setfontsize\normalsize\@xiipt{14.5}}
    \def\small{\@setfontsize\small\@xipt{13.6}}
    \let\footnotesize=\small
    \def\large{\@setfontsize\large\@xivpt{18}}
    \def\Large{\@setfontsize\Large\@xviipt{22}}
    \skip\@mpfootins = 18\p@ \@plus 2\p@
    \normalsize
} \@ifundefined{square}{}{} \makeatother

\def\file#1{\texttt{#1}}

\begin{document}
\let\WriteBookmarks\relax
\def\floatpagepagefraction{1}
\def\textpagefraction{.001}
\let\printorcid\relax 

\shorttitle{}    
\shortauthors{}  
\title [mode = title]{Error estimates of $hp$-finite element method for elliptic optimal control problems with Robin boundary}  

\tnotemark[1] 



\author[1]{Xingyuan Lin}

\credit{}


\author[1]{Xiuxiu Lin}

\author[1]{Xuesong Chen}
\cormark[1]
\ead{chenxs@gdut.edu.cn}

\credit{}

\affiliation[1]{organization={School of Mathematics and Statistics, Guangdong University of Technology},
            city={Guangzhou},
            postcode={510520}, 
            country={P.R. China}}

\cortext[1]{Corresponding author}



\begin{abstract}
  A priori and a posteriori error analysis of $hp$ finite element method for elliptic control problem with  Robin boundary condition and boundary observation are presented. are presented. Through the Clément-type approach and the construction of an auxiliary system, we derived a priori error estimates for the elliptic optimal control problem. Residual-based a posteriori error estimates are derived based on the well-known Scott-Zhang-type quasi-interpolation and coupled state-control approximations, thus establishing an a posteriori error estimator for the $hp$ finite element method. The numerical example demonstrates the accuracy of error estimation for the elliptic optimal control problems with Robin boundary.
\end{abstract}


\begin{highlights}
\item A priori and a posteriori error estimates for the $hp$-finite element method applied to elliptic optimal control with Robin boundary conditions and boundary observation.
\item Numerical experiments confirm the validity of the error estimator.
\end{highlights}



\maketitle

\section{Introduction}
In recent years, the study of optimal control problems has garnered attention across diverse fields and found extensive applications. For example, Grass et al. \cite{grass2019optimal} proposed a method integrating data-driven prediction with optimal control theory to efficiently manage and schedule fishery resources under coastal catch constraints. Zheng et al. \cite{zheng2024data} constructed an optimal control model to integrate submerged arc furnace systems with renewable energy sources (particularly photovoltaics and wind energy) and modern intelligent energy terminals, achieving dual benefits of economic efficiency and sustainable management. Amornraksa et al. \cite{amornraksa2024aopc} proposed a novel control strategy by incorporating multiple manipulated variables and biomass, enabling observation and control of wastewater concentration.

The finite element method has been widely employed in solving and analyzing error estimates for optimal control problems. Research in this area has evolved from regional control studies \cite{chen2011legendre,allendes2019adaptive,zhang2021conjugate} toward the investigation of boundary control \cite{wang2024weak,liu2024virtual,su2023numerical}. However, to the best of our knowledge, there remains limited research on error estimation for boundary state expectations under Robin boundary constraints. This paper aims to contribute to enriching this research area.

In optimal control problems, error analyses for different boundary conditions exhibit distinct characteristics. Dirichlet boundary conditions, being the most common type, have been extensively studied by numerous scholars \cite{langer2024adaptive,huang2014error,lin2023priori}. The literature \cite{leng2023} proposes an interior penalty hybridized discontinuous Galerkin (IP-HDG) method for elliptic problems, proves optimal and almost optimal approximation properties of the errors in the $W^{1,\infty}$ and $L^{\infty}$ norms, and applies the theoretical results to optimal control problems with pointwise state constraints. For Neumann boundary conditions, Liu et al. \cite{liu2024virtual} derived a priori error estimates in the $H^1$ and $L^2$ norms using a virtual element discretization scheme, first- and second-order optimality conditions, and auxiliary problems. Similar works include Wang et al. \cite{wang2024weak}, which obtained $L^2$ error estimates by enhancing the regularity of $u$ and the Lagrange multiplier $\lambda$. Mordukhovich \cite{mordukhovich2004neumann} analyzed hyperbolic optimal control problems with pointwise constraints on controls and states, establishing new necessary optimality conditions through the study of nonlinear multidimensional wave terms. Liu \cite{liu2016leapfrog} addresses parabolic PDE optimal control problems with the hybrid case of Robin boundary conditions, proposing a second-order finite difference scheme and rigorously proving its stability and convergence. Meanwhile, Ciaramella et al. \cite{ciaramella2022convergence} focused on applying the optimized Schwarz waveform relaxation method to periodic parabolic optimal control problems, adopting a semi-discrete time formulation and analyzing asymptotic behaviors as temporal discretization converges to zero.

When using the finite element method to study optimal control problems, solution accuracy can be enhanced by adjusting element size ($h$-version), increasing element order ($p$-version), or combining both ($hp$-version). The $h$-version finite element method, which refines element size, represents a traditional approach and is often integrated with adaptive mesh refinement strategies. For instance, Guo \cite{guo2023convergence} provided $L^2$ error analysis for sufficiently fine initial meshes and proved the convergence of adaptive finite element methods for nonmonotone elliptic problems; Winterscheidt \cite{winterscheidt2022improved} improved continuous $p$-order interpolation approximation functions to yield "semi-orthogonal" approximation functions, significantly reducing the condition number when generating stiffness and mass matrices for the Laplace equation. The $hp$-version method offers greater flexibility in mesh partitioning: the solution function tends to adjust element size $h$ in smooth regions and element order $p$ in nonsmooth regions. This characteristic endows the method with a natural advantage in improving adaptive finite element accuracy via a posteriori error estimators. Gong et al. \cite{gong2011posteriori} investigated a posteriori error estimates for the $hp$ finite element method and developed various residual-based weighted a posteriori error estimators. Chen and Lin \cite{chen2011posteriori} introduced new Clément-type and Scott–Zhang-type interpolation operators and studied a posteriori error estimates based on these new quasi-interpolation operators. Wachsmuth \cite{wachsmuth2016exponential} utilized $hp$-finite element discretization, analyzed solution regularity to demonstrate membership in appropriately chosen weighted Sobolev spaces, and established exponential convergence rates.

This paper investigates elliptic optimal control problems with Robin boundary conditions and boundary observations using the $hp$-finite element method. We derive both a priori and a posteriori error estimates, thus tackling the insufficiently studied challenges posed by Robin boundary constraints in control problems. The proposed residual-based estimators unify reliability and efficiency, facilitating adaptive refinement in practical implementations. These results expand the applicability of $hp$-FEM to boundary-driven optimization, providing theoretical insights for complex PDE-constrained scenarios.

The subsequent sections of this paper are structured as follows: section \ref{sec2} presents the problem formulation and preliminary results. Section \ref{sec3} derives a priori error estimates by leveraging discrete systems and continuous optimality conditions. Section \ref{sec4} constructs and analyzes the a posteriori error estimator, accompanied by proofs of its reliability and efficiency. Numerical experiments in section \ref{sec5} to verify the effectiveness of the estimator. Conclusion in section \ref{sec6}.

\section{$hp$ finite element approximation and preliminary results}\label{sec2}
In this paper, we use the $hp$-version finite element method for elliptic optimal control problems as follows:
\begin{equation}
\min_{u \in U_{ad}} J(y, u) = \frac{\lambda_{\Omega}}{2} \|y - y_{\Omega}\|_{L^{2}(\Omega)}^{2} + \frac{\lambda_{\Gamma}}{2} \|y - y_{\Gamma}\|_{L^{2}(\Gamma)}^{2} + \frac{\lambda}{2} \|u\|_{L^{2}(\Omega)}^{2} \label{eq1},
\end{equation}
subject to the elliptic partial differential equation
\begin{equation}
-\Delta y=\beta u\quad \text{in} \quad\Omega,\label{eq2}
\end{equation}
\begin{equation}
\partial_ny+\alpha y=0\quad \text{on} \quad\Gamma,\label{eq3}
\end{equation}
and the control constraints
\begin{equation}
U_{a}=\left\{u\in L^2(\Omega):u_{a}\leq u\quad a.e.\quad in\quad \Omega\right\}.\label{eq4}
\end{equation} 
Let $\Gamma$ be the boundary of domain $\Omega$, and $n$ denote the outward normal vector of $\Omega$, $\lambda_{\Omega}$ and $\lambda_{\Gamma}$ are non-negative constants, $\lambda$ are strictly positive constants.

Furthermore, $\beta\in L^{\infty}(\Omega)$ and $\alpha\in L^{\infty}(\Gamma)$. Here, $\alpha\geq M\textgreater 0$ almost everywhere on $\Gamma$, with $M$ being a constant. The target functions $y_{\Omega}$ and $y_{\Gamma}$ belong to $L^2(\Omega)$ and $L^2(\Gamma)$, respectively. Throughout the paper, the generic constants that appear in the theorems and lemmas are denoted by $C$.

Introduce the standard notation $W^{m,q}(\Omega)$ for Sobolev spaces on $\Omega$ equipped with the norm $\|\cdot\|_{W^{m,q}(\Omega)}$, the Sobolev space notation for the boundary $\Gamma$ is analogously defined. Specifically, $W^{m,q}(\Gamma)$ denotes the Sobolev space on $\Gamma$ equipped with the norm $\|\cdot\|_{W^{m,q}(\Gamma)}$. We denote $W^{m,q}(\Omega)$ by $H^{m}(\Omega)$ and $W^{m,q}(\Gamma)$ by $H^{m}(\Gamma)$.

To analyze the discrete problem in this section, the domain $\Omega\subset\mathbb{R}^2$ is defined as a convex
domain. We set $U=L^2(\Omega)$ as the control space and $Y=H^1(\Omega)$ as the state space. By applying Green's formula to (\ref{eq2}), we derive the weak formulation, and then substitute (\ref{eq3}) into it to obtain
\begin{equation}
\int_\Omega \nabla y\nabla v+\int_{\Gamma}\alpha yv=\int_{\Omega}\beta uv,\nonumber
\end{equation}
where $v\in Y$. We define
\begin{equation}
\begin{aligned}
&a(y,v)=\int_\Omega(\nabla y)\nabla v+\int_{\Gamma}(\alpha yv),\quad\forall y,v\in Y,\\
&\text{and}\\
&(u,v)_{A}=\int_{A} uv,\quad\forall v\in Y,\quad A= \Omega\text{ or }\Gamma.\end{aligned}\nonumber
\end{equation}
These equalities lead to the weak formulation of the control problem, which allows us to find $(y,u)\in Y\times U$,
\begin{equation}
\begin{aligned}
&\min_{u\in U_{ad}}J(y,u),\\
&a(y,v)=(\beta u,v)_{\Omega}.
\end{aligned}
\label{change_problem}
\end{equation}
Based on \cite{troeltzsch2010optimal}, the existence and uniqueness of the solution to the control problem \eqref{change_problem} are established. There exist constants $C,c\textgreater 0$ such that $\forall v,y\in Y$,
\begin{equation}
a(v,v)\geq c\|v\|_{H^1(\Omega)}^2,\quad a(y,v)\leq C\|y\|_{H^1(\Omega)}\|v\|_{H^1(\Omega)}.\nonumber
\end{equation}

\begin{Lemma}\label{yijiewenzi}
	(see \cite{troeltzsch2010optimal}) The control problem \eqref{change_problem} has a unique solution $(y,u)\in Y\times U$ if and only if there exists $z\in Y$ such that the triplet $(u,y,z)$ satisfies the following optimality conditions:
	
	\noindent
	\begin{minipage}{\linewidth}
		\begin{flalign}
		&\text{(i)} \quad a(y,v)=(\beta u,v)_{\Omega}, && v\in Y,&\nonumber\\
		&\text{(ii)} \quad a(q,z)=(\lambda_{\Omega}(y-y_{\Omega}),q)_{\Omega}+(\lambda_{\Gamma}(y-y_{\Gamma}),q)_{\Gamma}, && q\in Y,&\label{yijiexitong}\\
		&\text{(iii)} \quad (\beta z+\lambda u,w-u)\geq 0, && w \in U_{ad}.&\nonumber
		\end{flalign}
	\end{minipage}
\end{Lemma}

We discretize the problem using the $hp$-FEM. For further details, refer to \cite{melenk2005hp}. The domain $\Omega$ is partitioned into non-overlapping elements $\tau$, with the collection of these elements denoted by $\mathcal{T}$. Each element $\tau$ is associated with an affine mapping $F_{\tau}:\hat{\tau}\to\tau$, where $\hat{\tau}$ represents the reference element—specifically, $\hat{\tau}=(0,1)^2$. We define $h_{\tau}=\text{diam }\tau$ (the diameter of element $\tau$). The set of all edges of $\tau$ is denoted by $\mathcal{E}(\tau)$, and $\mathcal{E}(\mathcal{T})$ represents the set of all edges in the triangulation $\mathcal{T}$. Additionally, the set of all vertices of $\tau$ is denoted by $\mathcal{N}(\tau)$, while $\mathcal{N}(\mathcal{T})$ represents the set of all vertices in the triangulation $\mathcal{T}$.\\
We assume that the triangulation is $\gamma\operatorname{-shape}$ regular,
\begin{equation}
h_\tau{}^{-1}\|F_\tau^{\prime}\|+h_\tau\|(F_\tau^{\prime})^{-1}\|\leq\gamma.\nonumber
\end{equation}
For each element $\tau$, there corresponds a polynomial degree $p_{\tau}\in\mathbb{N}$. We collect all $p_{\tau}$ into a vector  $\mathbf{p}=(p_\tau)_{\tau\in\mathcal{T}}$. We can define the spaces $U^{\mathbf{p}}(\mathcal{T},\Omega)$ and $ S^{\mathbf{p}}(\mathcal{T},\Omega)$ as follows:
\begin{equation}
\begin{aligned}
U^{\mathbf{p}}(\mathcal{T},\Omega)&:=\{u\in L^{2}(\Omega): u\lvert_{\tau} \circ F_{\tau}\in \Pi_{p}(\hat{\tau})\},\\
S^{\mathbf{p}}(\mathcal{T},\Omega)&:=\{v\in H^{1}(\Omega): v\lvert_{\tau} \circ F_{\tau}\in \Pi_{p}(\hat{\tau})\}.
\end{aligned}
\nonumber
\end{equation}
We set
\begin{equation}
\Pi_k(\hat{\tau}) = \begin{cases}
P_k := \operatorname{span}\{x^i y^j : 0 \leq i + j \leq k\}, & \operatorname{if}\text{ }\hat{\tau} = T, \\
Q_k := \operatorname{span}\{x^i y^j : 0 \leq i, j \leq k\}, & \operatorname{if}\text{ }\hat{\tau} = S,
\end{cases} \nonumber
\end{equation}
here, $T$ represents a triangle and $S$ represents a rectangle.

We will write $S^\mathbf{p}(\mathcal{T})$ if the degree vector $\mathbf{p}$ satisfes $p_{\tau
}=p$ for all $\tau \in \mathcal{T}$. Moreover, we need to assume that for any two adjacent elements $\tau$ and $\tau^{\prime}$, the polynomial degree vector $\mathbf{p}$ satisfies\\
\begin{equation}
\gamma^{-1}p_{\tau} \leq p_{\tau^{\prime}} \leq \gamma p_{\tau} \quad \tau,\tau^{\prime}\in\mathcal{T}\mathrm{~with~}\overline{\tau}\cap\overline{\tau}^{\prime}\neq\emptyset.
\label{degree_constraint}
\end{equation}
We obtain the $hp$-FEM state space $Y^{hp}=S^{\mathbf{p}}(\mathcal{T},\Omega)\cap Y$ and the $hp$-FEM control space $U^{hp}=U^{p}(\mathcal{T},\Omega)\cap U_{ad}$ associated with the triangulation $\mathcal{T}$.

\noindent Consequently, the $hp$ finite element discretization of the original problem can be rewritten as
\begin{equation}
\min_{u_{hp}\in U^{hp}}J(y_{hp},u_{hp})=\frac{\lambda_{\Omega}}{2}\|y_{hp}-y_{\Omega}\|_{L^{2}(\Omega)}^{2}+\frac{\lambda_{\Gamma}}{2}\|y_{hp}-y_{\Gamma}\|_{L^{2}(\Gamma)}^{2}+\frac{\lambda}{2}\|u_{hp}\|_{L^{2}(\Omega)}^{2},\label{2.2}
\end{equation}
\begin{equation}
a(y_{hp},v_{hp})=(\beta u_{hp},v_{hp})_{\Omega}, \quad v_{hp}\in Y^{hp}.\label{2.3}
\end{equation}
Since $J(y_{hp},u_{hp})$ is a strongly convex function, the solution to constrained optimization problem $\eqref{2.2}-\eqref{2.3}$ exists and is unique.
\begin{Lemma}\label{lisanwenzi}
	For the discrete control problem ($\ref{2.2}$)-($\ref{2.3}$), there exists a unique solution $(y_{hp},u_{hp})\in Y^{hp}\times U^{hp}$ if and only if there exists $z_{hp}\in Y^{hp}$ such that the triplet $(y_{hp},u_{hp},z_{hp})$ satisfies the following equations:
	
	\noindent
	\begin{minipage}{\linewidth}
		\begin{flalign}
		&\text{(i)}a(y_{hp},v_{hp})=(\beta u_{hp},v_{hp})_{\Omega}, && v_{hp}\in Y^{hp},&\nonumber\\
		&\text{(ii)} a(q_{hp},z_{hp})=(\lambda_{\Omega}(y_{hp}-y_{\Omega}),q_{hp})_{\Omega}+(\lambda_{\Gamma}(y_{hp}-y_{\Gamma}),q_{hp})_{\Gamma}, && q_{hp}\in Y^{hp},&\nonumber\\
		&\text{(iii)} (\beta z_{hp}+\lambda u_{hp},w_{hp}-u_{hp})\geq 0, && w_{hp} \in U^{hp}.&\label{lisanxitong}
		\end{flalign}
	\end{minipage}
\end{Lemma}

\section{A priori error estimate}\label{sec3}
In this section, we present a priori error estimate. To derive the estimate, it is essential to introduce the Clément-type interpolation operator, construct an auxiliary system, and prove several theoretical lemmas. \\
The error bounds are derived by employing the standard interpolation operator based on the reference domain, techniques from the $hp$-finite element method, and the mapping properties of state and control variables.\\
We define the projection operator $\pi_{p}^{h}$ as follows: $\forall v\in Y$, find $\pi_{p}^{h}v\in S^{\mathbf{p}}(\mathcal{T},\Omega)$ such that
\begin{equation}
a(\pi_{p}^{h}v-v,v_{hp})=0,\quad\forall\quad\upsilon_{hp}\in S^{\mathbf{p}}(\mathcal{T},\Omega).
\nonumber
\end{equation}
\begin{Lemma}\label{Byinli}
	(see \cite{babuska1987hp})Let $\pi_p^h:Y\cap H^s(\Omega)\longrightarrow S^{\mathbf{p}}(\mathcal{T},\Omega)$ such that for any $0\leq t\leq s$, $\forall u\in Y\cap H^s$
	\begin{equation}
	\|u-\pi_{p}^{h}u\|_{H^t(\Omega)}\leq C\frac{h^{\mu-t}}{p^{s-t}}\|u\|_{H^{s}(\Omega)},
	\label{Celemnt_type1}
	\end{equation}
	where $\mu=min\{p+1,s\}$. 
\end{Lemma}
\noindent We introduce auxiliary variables $y_{hp}(u),z_{hp}(u)\in Y^{hp}$, defined as the solutions to the following equations:
\begin{equation}
\begin{aligned}
(i)&a(y_{hp}(u),v_{hp})=(\beta u,v_{hp})_{\Omega}, \quad v_{hp}\in Y^{hp},\\
(ii)&a(q_{hp},z_{hp}(u))=(\lambda_{\Omega}(y-y_{\Omega}),q_{hp})_{\Omega}+(\lambda_{\Gamma}(y-y_{\Gamma}),q_{hp})_{\Gamma},  \quad q_{hp}\in Y^{hp}.\\
\end{aligned}
\label{Xianyanfuzhuxitong1}
\end{equation}
\begin{Lemma}\label{fisrt_part_of_Xianyan}
	Let $(u,y,z)$ be the optimal solution of the problem \eqref{yijiexitong}, $(y_{hp}(u),z_{hp}(u))$ be the solution of the auxiliary system \eqref{Xianyanfuzhuxitong1}, then there holds the following estimate:
	\begin{equation}
	\|y-y_{hp}(u)\|_{H^{1}(\Omega)}+\|z-z_{hp}(u)\|_{H^1(\Omega)}\leq C\frac{h^{\mu-1}}{p^{m-1}}\left(\|y\|_{H^{m}(\Omega)}+\|z\|_{H^{m}(\Omega)}\right),
	\label{first_part_of_Xianyaneq}
	\end{equation}
	where $\mu=min\{p+1,m\}$.
\end{Lemma}
\noindent Proof. Combining the optimality conditions \eqref{yijiexitong} and auxiliary system \eqref{Xianyanfuzhuxitong1}, we have
\begin{equation}
a(y-y_{hp}(u),v_{hp})=0 \quad \forall v_{hp}\in Y^{hp},
\label{fisrt_part_of_Xianyan_y1}
\end{equation}
\begin{equation}
a(q_{hp},z-z_{hp}(u))=0 \quad \forall q_{hp}\in Y^{hp}.
\label{fisrt_part_of_Xianyan_p1}
\end{equation}
From \eqref{fisrt_part_of_Xianyan_y1}, $\forall w_{hp}\in Y^{hp}$ satisfied 
\begin{equation}
\begin{aligned}
&c\|y-y_{hp}(u)\|_{H^1(\Omega)}^2\\
&\leq a(y-y_{hp}(u),y-w_{hp})+a(y-y_{hp}(u),w_{hp}-y_{hp}(u))\\
&=a(y-y_{hp}(u),y-w_{hp})\\
&\leq C\|y-y_{hp}(u)\|_{H^1(\Omega)}\|y-w_{hp}\|_{H^1(\Omega)}\\
&\leq C\|y-y_{hp}(u)\|_{H^1(\Omega)}\inf_{w_{hp}\in Y^{hp}}\|y-w_{hp}\|_{H^1(\Omega)}.
\nonumber
\end{aligned}
\end{equation}
Canceling the common terms of the above inequality, and using \eqref{Celemnt_type1}
\begin{equation}
\begin{aligned}
\|y-y_{hp}(u)\|_{H^1(\Omega)}&\leq C\inf_{w_{hp}\in Y^{hp}}\|y-w_{hp}\|_{H^1(\Omega)}\\
&\leq C\|y-\pi_{p}^{h}y\|_{H^1(\Omega)}\\
&\leq C\frac{h^{\mu-1}}{p^{m-1}}\|y\|_{H^{m}(\Omega)}.
\end{aligned}
\label{first_yHm}
\end{equation}
By analogous methodology with \eqref{fisrt_part_of_Xianyan_p1}, $\forall q_{hp}\in Y^{hp}$, we have
\begin{equation}
\begin{aligned}
c\|z-z_{hp}(u)\|_{H^1(\Omega)}^2&\leq a(z-q_{hp},z-z_{hp}(u))+a(q_{hp}-z_{hp},z-z_{hp}(u))\\
&\leq C\|z-z_{hp}(u)\|_{H^1(\Omega)}\inf_{q_{hp}\in Y^{hp}}\|z-q_{hp}\|_{H^1(\Omega)}
\nonumber
\end{aligned}
\end{equation}
From \eqref{Celemnt_type1}, we obtain
\begin{equation}
\begin{aligned}
\|z-z_{hp}(u)\|_{H^1(\Omega)}&\leq C\inf_{q_{hp}\in Y^{hp}}\|z-q_{hp}\|_{H^1(\Omega)}\\
&\leq C\|y-\pi_{p}^{h}z\|_{H^1(\Omega)}\\
&\leq C\frac{h^{\mu-1}}{p^{m-1}}\|z\|_{H^{m}(\Omega)}.
\end{aligned}
\label{first_pHm}
\end{equation}
Then the proof complete. $\hfill\square$
\begin{Lemma}\label{second_part_of_Xianyan}
	Let $(u_{hp},y_{hp},z_{hp})$ and $(y_{hp}(u),z_{hp}(u))$ be the optimal solution of the problem \eqref{lisanxitong} and \eqref{Xianyanfuzhuxitong1}, respectively. There holds the following estimate:
	\begin{equation}
	\begin{aligned}
	\|y_{hp}-y_{hp}(u)\|_{H^{1}(\Omega)}&+\|z_{hp}-z_{hp}(u)\|_{H^1(\Omega)}\\
	&\leq C\left(\|u_{hp}-u\|_{L^{2}(\Omega)}+\|y-y_{hp}\|_{L^2(\Omega)}\right).
	\end{aligned}
	\nonumber
	\end{equation}
\end{Lemma}
\noindent Proof.\quad Combining the discrete systems \eqref{lisanxitong} and auxiliary system \eqref{Xianyanfuzhuxitong1}, we have
\begin{equation}
a(y_{hp}-y_{hp}(u),v_{hp})=(\beta(u_{hp}-u),v_{hp})_{\Omega}\quad\forall v_{hp}\in Y^{hp},
\label{second_part_XY_y1}
\end{equation}
\begin{equation}
a(q_{hp},z_{hp}-z_{hp}(u))=(\lambda_{\Omega}(y_{hp}-y),q_{hp})_{\Omega}+(\lambda_{\Gamma}(y_{hp}-y),q_{hp})_{\Gamma}\quad\forall q_{hp}\in Y^{hp}.
\label{second_part_XY_p1}
\end{equation}
Letting $v_{hp}=y_{hp}-y_{hp}(u)$, from \eqref{second_part_XY_y1}, we get
\begin{equation}
\begin{aligned}
c\|y_{hp}-y_{hp}(u)\|_{H^1(\Omega)}^2&\leq a(y_{hp}-y_{hp}(u),y_{hp}-y_{hp}(u)) \\
&=(\beta(u_{hp}-u),y_{hp}-y_{hp}(u))_{\Omega}\\
&\leq C\|u_{hp}-u\|_{L^2(\Omega)}\|y_{hp}-y_{hp}(u)\|_{H^1(\Omega)},
\nonumber
\end{aligned}
\end{equation}
thus we derive
\begin{equation}
\|y_{hp}-y_{hp}(u)\|_{H^1(\Omega)}\leq C\|u_{hp}-u\|_{L^2(\Omega)}.\nonumber
\end{equation}
Let $q_{hp}=p_{hp}-z_{hp}(u)$, from \eqref{second_part_XY_p1}, we get
\begin{equation}
\begin{aligned}
&c\|z_{hp}-z_{hp}(u)\|_{H^1(\Omega)}^2\leq a(z_{hp}-z_{hp}(u),z_{hp}-z_{hp}(u)) \\
&=(\lambda_{\Omega}(y_{hp}-y),z_{hp}-z_{hp}(u))_{\Omega}+(\lambda_{\Gamma}(y_{hp}-y),z_{hp}-z_{hp}(u))_{\Gamma}\\
&\leq C\|y_{hp}-y\|_{L^2(\Omega)}\|z_{hp}-z_{hp}(u)\|_{L^2(\Omega)}\\
&+C\|y_{hp}-y\|_{L^2(\Gamma)}\|z_{hp}-z_{hp}(u)\|_{L^2(\Gamma)}\\
&\leq C\|y_{hp}-y\|_{L^2(\Omega)}\|z_{hp}-z_{hp}(u)\|_{H^1(\Omega)}.
\nonumber
\end{aligned}
\end{equation}
thus we have
\begin{equation}
\|z_{hp}-z_{hp}(u)\|_{H^1(\Omega)}\leq C\|y_{hp}-y\|_{L^2(\Omega)}.\nonumber
\end{equation}
Then the proof is complete. $\hfill\square$\\
Based on the above lemmas, we proceed to present the first estimate in Theorem 1.
\begin{Theorem}\label{Xianyan_u_estimate}
	Let $(u,y,z)$ and $(u_{hp},y_{hp},z_{hp})$ denote the optimal solutions of the problem \eqref{yijiexitong} and \eqref{lisanxitong}, respectively. And $(y_{hp}(u),z_{hp}(u))$ be the solution of the auxiliary system \eqref{Xianyanfuzhuxitong1}, then there holds the following estimate:
	\begin{equation}
	\nu\|u-u_{hp}\|_{L^2(\Omega)}^2+\|y-y_{hp}\|_{L^2(\Omega)}^2\leq C_{\nu}\|z-z_{hp}(u)\|_{H^1(\Omega)}^2+\|y-y_{hp}(u)\|_{H^1(\Omega)}^2,
	\nonumber
	\end{equation}
	where $\nu$ depends solely on the constant $\lambda$.
\end{Theorem}
\noindent Proof. Letting $v_{hp}=z_{hp}(u)-z_{hp}$ in \eqref{second_part_XY_y1},\quad $q_{hp}=y_{hp}-y_{hp}(u)$ in \eqref{second_part_XY_p1}, using trace theorem, we have
\begin{equation}
\begin{aligned}
&(\beta(u_{hp}-u),z_{hp}(u)-z_{hp})_{\Omega}\\
&=(\lambda_{\Omega}(y-y_{hp}),y_{hp}-y_{hp}(u))_{\Omega}+(\lambda_{\Gamma}(y-y_{hp}),y_{hp}-y_{hp}(u))_{\Gamma}\\
&=\frac{\lambda_{\Omega}}{2}\{-\|y-y_{hp}\|_{L^2(\Omega)}^2+\|y-y_{hp}(u)\|_{L^2(\Omega)}^2-\|y_{hp}-y_{hp}(u)\|_{L^2(\Omega)}^2\}\\
&+\frac{\lambda_{\Gamma}}{2}\{-\|y-y_{hp}\|_{L^2(\Gamma)}^2+\|y-y_{hp}(u)\|_{L^2(\Gamma)}^2-\|y_{hp}-y_{hp}(u)\|_{L^2(\Gamma)}^2\}\\
&\leq C\{\|y-y_{hp}(u)\|_{L^2(\Omega)}^2+\|y-y_{hp}(u)\|_{L^2(\Gamma)}^2-\|y-y_{hp}\|_{L^2(\Omega)}^2\}\\
&\leq C\{\|y-y_{hp}(u)\|_{H^1(\Omega)}^2-\|y-y_{hp}\|_{L^2(\Omega)}^2\}.
\end{aligned}
\label{y-y_hp_ofL2}
\end{equation}
Letting
\begin{equation}
j_{hp}=\left\{
\begin{array}{l}
cu_{hp},\quad u_{a}\leq c u_{hp}\\
c_{min}u_{hp},\quad c u_{hp}<u_{a}.
\end{array}\right.
\nonumber
\end{equation}
where $c=\frac{(u, \lambda u_{hp}+\beta z_{hp})}{(u_{hp}, \lambda u_{hp}+\beta z_{hp})},\quad c_{min}=\frac{u_{a}}{min\left\{u_{hp}\right\}}$, 
$v_{hp}=\frac{(u,\lambda u_{hp}+\beta z_{hp})}{(u_{hp},\lambda u_{hp}+\beta z_{hp})}u_{hp}\in U^{hp}$.
We thus obtain two cases,
\begin{align}
Case\quad I: u_{a}\leq cu_{hp}&\Rightarrow (u-j_{hp}, \beta z_{hp}+\lambda u_{hp})=0,\nonumber\\
Case\quad II: cu_{hp}<u_{a}&\Rightarrow (u-cu_{hp}+cu_{hp}-j_{hp}, \beta z_{hp}+\lambda u_{hp}),\nonumber\\
&=0+(c-c_{min})(u_{hp}, \beta z_{hp}+\lambda u_{hp})\geq 0,\nonumber
\end{align}
such that $(u-j_{hp},\beta z_{hp}+\lambda u_{hp})\geq0$, from (iii) of \eqref{yijiexitong}, (iii) of \eqref{lisanxitong}, \eqref{y-y_hp_ofL2}, and Young’s inequality,
\begin{equation}
\begin{aligned}
&\lambda\|u-u_{hp}\|_{L^2(\Omega)}^2=(u-u_{hp},\lambda u-\lambda u_{hp})\\
&=(u-u_{hp},\lambda u+\beta z-\beta z+\beta z_{hp}-\beta z_{hp}-\lambda u_{hp})\\
&=(u-u_{hp},\lambda u+\beta z)-(u-u_{hp},\beta z_{hp}+\lambda u_{hp})+(u-u_{hp},\beta z_{hp}-\beta z)\\
&\leq(u-u_{hp},\beta z_{hp}-\beta z)-(u-u_{hp},\beta z_{hp}+\lambda u_{hp})\\
&\leq(u-u_{hp},\beta z_{hp}-\beta z_{hp}(u)+\beta z_{hp}(u)-\beta z)-\left(u-j_{hp}+j_{hp}-u_{hp},\beta z_{hp}+\lambda u_{hp}\right)\\
&\leq(u-u_{hp},\beta z_{hp}-\beta z_{hp}(u))+(u-u_{hp},\beta z_{hp}(u)-\beta z)-\left(u-j_{hp},\beta z_{hp}+\lambda u_{hp}\right)+\left(u_{hp}-j_{hp},\beta z_{hp}+\lambda u_{hp}\right)\\
&\leq(u-u_{hp},\beta z_{hp}-\beta z_{hp}(u))+(u-u_{hp},\beta z_{hp}(u)-\beta z)-(u-v_{hp},\beta z_{hp}+\lambda u_{hp})\\
&\leq (u-u_{hp},\beta z_{hp}-\beta z_{hp}(u))+(u-u_{hp},\beta z_{hp}(u)-\beta z)\\
&\leq \|u_{hp}-u\|_{L^2(\Omega)}\|\beta(z-z_{hp}(u))\|_{L^2(\Omega)}+C\{\|y-y_{hp}(u)\|_{L^1(\Omega)}^2-\|y-y_{hp}\|_{L^2(\Omega)}^2\}\\
&\leq \frac{1}{2}\|u_{hp}-u\|_{L^2(\Omega)}^2+\frac{\beta^2}{2}\|(z-z_{hp}(u))\|_{L^2(\Omega)}^2+C\{\|y-y_{hp}(u)\|_{L^1(\Omega)}^2-\|y-y_{hp}\|_{L^2(\Omega)}^2\},
\end{aligned}
\nonumber
\end{equation} 
letting $\nu=\frac{2\lambda-1}{2C},C_{\nu}=\frac{\beta^2\nu}{2\lambda-1}$, we get
\begin{equation}
\nu\|u-u_{hp}\|_{L^2(\Omega)}^2+\|y-y_{hp}\|_{L^2(\Omega)}^2\leq C_{\nu}\|z-z_{hp}(u)\|_{H^1(\Omega)}^2+\|y-y_{hp}(u)\|_{H^1(\Omega)}^2.\nonumber\hfill\square
\end{equation}
\begin{Theorem}\label{xianyanwucha}
	Let $(y,u,z)$ and ($y_{hp},u_{hp},z_{hp}$) be the solution of $\eqref{yijiexitong}$ and $\eqref{lisanxitong}$ respectively. Assuming that $(y,z)\in H^{m}(\Omega)\times H^{m}(\Omega)(m\geq 1)$, we obtain a priori error estimate as follows
	\begin{equation}
	\begin{aligned}
	\|u-u_{hp}\|_{L^2(\Omega)}+\|y-y_{hp}\|_{H^1(\Omega)}&+\|z-z_{hp}\|_{H^1(\Omega)}\\&\leq C\frac{h^{\mu-1}}{p^{m-1}}\{\|y\|_{H^{m}(\Omega)}+\|z\|_{H^{m}(\Omega)}\},
	\end{aligned}
	\nonumber
	\end{equation}
	where $\mu=min\{p+1,m\}$.
\end{Theorem}
\noindent Proof. From Lemma \ref{second_part_of_Xianyan} and Theorem \ref{Xianyan_u_estimate}, \eqref{first_yHm}, \eqref{first_pHm}, we get
\begin{equation}
\begin{aligned}
&\|y_{hp}-y_{hp}(u)\|_{H^1(\Omega)}+\|z_{hp}-z_{hp}(u)\|_{H^1(\Omega)}\\
&\leq C\{\|u_{hp}-u\|_{L^2(\Omega)}+\|y-y_{hp}\|_{L^2(\Omega)}\}\\
&\leq C\frac{h^{\mu-1}}{p^{m-1}}\sqrt{\|z-z_{hp}(u)\|_{H^1(\Omega)}^2+\|y-y_{hp}(u)\|_{H^1(\Omega)}^2}\\
&\leq C\frac{h^{\mu-1}}{p^{m-1}}\sqrt{\|z\|_{H^{m}(\Omega)}^2+\|y\|_{H^{m}(\Omega)}^2}\\
&\leq C\frac{h^{\mu-1}}{p^{m-1}}\{\|z\|_{H^{m}(\Omega)}+\|y\|_{H^{m}(\Omega)}\}.
\end{aligned}
\label{Xianyan_fullpart}
\end{equation}
From \eqref{first_part_of_Xianyaneq} and \eqref{Xianyan_fullpart}, we have
\begin{equation}
\begin{aligned}
&\|u-u_{hp}\|_{L^2(\Omega)}+\|y-y_{hp}\|_{H^1(\Omega)}+\|z-z_{hp}\|_{H^1(\Omega)}\\
&\leq \|u-u_{hp}\|_{L^2(\Omega)}+\|y-y_{hp}(u)\|_{H^{1}(\Omega)}+\|z-z_{hp}(u)\|_{H^1(\Omega)}\\
&+\|y_{hp}-y_{hp}(u)\|_{H^{1}(\Omega)}+\|z_{hp}-z_{hp}(u)\|_{H^1(\Omega)}\\
&\leq C\frac{h^{\mu-1}}{p^{m-1}}(\|y\|_{H^{m}(\Omega)}+\|z\|_{H^{m}(\Omega)}). 
\nonumber
\end{aligned}
\end{equation}
The proof is complete.$\hfill\square$\\
The a priori error analysis developed in this section establishes theoretical convergence rates for the $hp$-FEM discretization of the Robin boundary control problem. While these results provide valuable insight into the asymptotic behavior of the method, their practical utility is limited by their dependence on the unknown solution regularity and global discretization parameters. To bridge this gap, we now turn to a posteriori error estimation.

\section{A posterior error estimate and algorithm}\label{sec4}
\subsection{A posterior error estimate}
In this section, we derive upper bounds for the a posteriori error. By leveraging the Scott–Zhang-type quasi-interpolation commonly used in $hp$-FEM, we draw upon the continuity conditions from \cite{melenk2005hp} to establish the approximation results.\\
The boundary edge set $\mathcal{E}(\mathcal{T})$ of the triangulation $\mathcal{T}$ satisfies $\mathcal{E}_{0}(\mathcal{T}) \subset \mathcal{E}(\mathcal{T})$, i.e.,
\begin{equation}
\quad \mathcal{E}_{0}(\mathcal{T}) \subset \mathcal{E}(\mathcal{T}) \quad 
\text{and} \quad b \in \partial\Omega,\ \forall b \in \mathcal{E}_{0}(\mathcal{T}).
\nonumber
\end{equation}
Define for $q \in (1,\infty)$:
\begin{equation}
\begin{aligned}
&W_{\mathcal{E}_{0}(\mathcal{T}),\mathbf{p}}^{1,q} :=\{ u \in W^{1,q}(\Omega) :u\lvert_b \circ F_b \in \mathcal{P}_{p_b},\forall b, b' \in \mathcal{E}_{0}(\mathcal{T}),V \in \Lambda(b) \cap \Lambda(b'),\\&\lim_{\substack{x \to V \\ x \in b}} u(x) = \lim_{\substack{x \to V \\ x \in b'}} u(x)\},
\end{aligned}
\label{continuity condition}
\end{equation}
where $\mathcal{P}_{p_b}$ denotes the space of polynomials of degree at most $p_b$ defined on the reference element associated with the boundary edge $b,\Lambda(b)$ is the set of vertices of the boundary edge $b$, and $V$ is a common vertex shared by edges $b$ and $b'$.\\
We introduce the notion of patch $\omega_V$ associated with the node $V\in \mathcal{N}(\mathcal{T})$ by
\begin{equation}
\omega_V:=\{x\in\Omega\lvert x\in\overline{\tau}\quad\text{ for some }\tau\mathrm{~with~}V\in\overline{\tau}\}^\circ,\nonumber
\end{equation}
where $A^\circ$ denotes the interior of the set $A$. Patches of order $j\in\mathbb{N}$ associated with an element $\tau\in \mathcal{T}$ or an edge $e\in\mathcal{E}(\mathcal{T})$ are defined thus:
\begin{equation}
\begin{aligned}
\omega_e^1:&=\bigcup_{V\in\mathcal{N}(e)}\omega_V,\quad\omega_e^{j+1}:=\bigcup_{V\in\mathcal{N}(\mathcal{T}):V\in\overline{\omega_e^j}}\omega_V,\quad j=1,2,\ldots,\\
\omega_\tau^1:&=\bigcup_{V\in\mathcal{N}(\tau)}\omega_V,\quad\omega_\tau^{j+1}:=\bigcup_{V\in\mathcal{N}(\mathcal{T}):V\in\overline{\omega_\tau^j}}\omega_V,\quad j=1,2,\ldots.
\end{aligned}
\nonumber
\end{equation}
\begin{Lemma}\label{Scott1wenzi}
	(Scott-Zhang-type quasi interpolation\cite{melenk2005hp}). Let $\mathcal{T}$ be a $\gamma$-shape regular triangulation of a domain $\Omega\subset\mathbb{R}^2$. Let $\mathbf{p}$ be a polynomial degree distribution satisfying \eqref{degree_constraint}. Let $\mathcal{E}_{0}(\mathcal{T})\subset\mathcal{E}(\mathcal{T})$ be a collection of boundary edges and $q\neq 2$ assume additionally that $\|p_{\tau}-p_{\tau^{\prime}}\|\leq \gamma,\forall\tau,\tau^{\prime }\text{s.t.} \overline {\tau}\cap \overline {\tau^{\prime }}\cap \overline {b}\neq \emptyset$ for some $b\in\mathcal{E}_{0}(\mathcal{T})$. Let the continuity condition \eqref{continuity condition} be satisfied. Then there exists a linear operator $I_{1}^{hp}: H^{1}(\Omega) \to S^{\mathbf{p}}(\mathcal{T},\Omega)$ such that
	\begin{equation}
	(I_{1}^{hp}u)\lvert_b=u\lvert_b,\quad\forall b\in\mathcal{E}_{0}(\mathcal{T}).
	\nonumber
	\end{equation}
	Furthermore, there exists a constant $C>0$ depending only on $\gamma$ and $q$ such that for all elements $\tau\in\mathcal{T}$ and all edges $e\in \mathcal{E}(\mathcal{T})$,
	\begin{align}
	&\|u-I_{1}^{hp}u\|_{L^2(\tau)}+\frac{h_\tau}{p_\tau}\|\nabla(u-I_{1}^{hp}u)\|_{L^2(\tau)}\leq C\frac{h_\tau}{p_\tau}\|\nabla u\|_{L^2(\omega_\tau^4)},\nonumber\\
	&\|u-I_{1}^{hp}u\|_{L^2(e)}\leq C\left(\frac{h_\tau}{p_\tau}\right)^{1/2}\|\nabla u\|_{L^2(\omega_e^4)}.\nonumber
	\end{align}
\end{Lemma}
\begin{Lemma}\label{Ihp1wenzi}
	(see \cite{chen2011posteriori}). There exists a constant C>0 independent of $v$, $h_{\tau}$ and $p_{\tau}$ and a mapping $\pi_{p_{\tau}}^{h_{\tau}}:H^1(\tau)\to \mathscr{P}$ such that $v\in H^1(\tau)$,$\tau\in\mathcal{T}$ the following inequality is valid
	\begin{equation}
	\|v-\pi_{p_{\tau}}^{h_{\tau}}v\|_{L^2(\tau)}\leq C\frac{h_{\tau}}{p_{\tau}}\|v\|_{H^1(\tau)},
	\nonumber
	\end{equation}
	where we will write $v\in\mathscr{P}$ if the following satisfied: $v\lvert_{\tau}\circ F_{\tau}\in P_{p_{\tau}}(\hat{\tau})$ is a triangle; $v\lvert_{\tau}\circ F_{\tau}\in Q_{p_{\tau}}(\hat{\tau})$ is a parallelogram.\\
	Let $H^*(\Omega,\mathcal{T})=\{v:v\lvert_{\tau}\in H^1(\tau),\forall\tau\in\mathcal{T}\}$ and then we can define the mapping that is useful in the estimate of the
	control, i.e., there exists a mapping $I_U^{hp}:H^*(\Omega,\mathcal{T})\to U^{\mathbf{p}}(\mathcal{T})$ such that \\
	\begin{equation}
	I_U^{hp}v\lvert_{\tau}=\pi_{p_{\tau}}^{h_{\tau}}(v\lvert_{\tau}),\quad\forall\tau\in\mathcal{T}.
	\nonumber
	\end{equation}
\end{Lemma}
\noindent We investigate the posteriori error for problem \eqref{change_problem}. We will introduce $y(u_{hp})$ and $z(u_{hp})$, defined by the auxiliary system:
\begin{equation}
\begin{aligned}
(i)\quad a(y(u_{hp}), v) &= (\beta u_{hp}, v)_{\Omega} \quad \forall v \in Y, \\
(ii)\quad a(q, z(u_{hp})) &= (\lambda_{\Omega}(y(u_{hp}) - y_{\Omega}), q)_{\Omega} + (\lambda_{\Gamma}(y(u_{hp}) - y_{\Gamma}), q)_{\Gamma} \quad \forall q \in Y.
\end{aligned}
\label{fuzhuxitong}
\end{equation}
\begin{Definition}\label{definition1xianyan}
	We define the notations:
	\begin{flalign}
	\eta^2 &= \sum_{i=1}^{7} \eta_{i}^2, & \nonumber \\
	\eta_{1}^2 &=\sum\limits_{\tau\in \mathcal{T}}\frac{h_\tau^2}{p_\tau^2}\|\beta u_{hp}+\Delta y_{hp}\|_{L^2(\Omega)}^2, & \nonumber \\
	\eta_{2}^2 &=\sum\limits_{e\in \mathcal{E}(\mathcal{T})\backslash\mathcal{E}_{0}(\mathcal{T})}\frac{h_e}{p_e}\|\nabla y_{hp}\cdot n_e\|_{L^2(e)}^2, & \nonumber \\
	\eta_{3}^2 &=\sum\limits_{e\in \mathcal{E}_{0}(\mathcal{T})}\frac{h_e}{p_e}\|\alpha y_{hp}+[\nabla y_{hp}\cdot n_e]\|_{L^2(e)}^2, & \nonumber \\
	\eta_{4}^2 &=\sum\limits_{\tau\in \mathcal{T}}\frac{h_\tau^2}{p_\tau^2}\|\lambda_{\Omega}(y_{hp}-y_{\Omega})+\Delta z_{hp}\|_{L^2(\tau)}^2, & \nonumber \\
	\eta_{5}^2 &=\sum\limits_{e\in \mathcal{E}(\mathcal{T})\backslash\mathcal{E}_{0}(\mathcal{T})}\frac{h_e}{p_e}\|[\nabla z_{hp}\cdot n_e]\|_{L^2(e)}^2, & \nonumber \\
	\eta_{6}^2 &=\sum\limits_{e\in \mathcal{E}_{0}(\mathcal{T})}\frac{h_e}{p_e}\|\lambda_{\Gamma}(y_{hp}-y_{\Gamma})-\alpha z_{hp}-[\nabla z_{hp}\cdot n_e]\|_{L^2(e)}^2, & \nonumber \\
	\eta_{7}^2 &=\sum\limits_{\tau\in \mathcal{T}}\frac{h_{\tau}^2}{p_{\tau}^2}\|\nabla(\lambda u_{hp}+\beta z_{hp})\|_{L^2(\tau)}^2, & \nonumber
	\end{flalign}
\end{Definition}
\noindent where $n_e$ is the unit outer normal vector of edge $e$.
From (iii) of \eqref{yijiexitong},
\begin{equation}
(\lambda u, u - u_{hp}) \leq -(\beta z, u - u_{hp}). 
\label{lambda_u_xiaoyu1}
\end{equation}
While (iii) of \eqref{lisanxitong} states the positivity:
\begin{equation}
(\beta z_{hp} + \lambda u_{hp}, w_{hp} - u_{hp}) \geq 0, \quad \forall w_{hp} \in U^{hp}.
\label{lambda_u_xiaoyu2}
\end{equation}
Following the error estimation techniques in \cite{chen2011posteriori}, we introduce the projection operator $\Pi^{hp}: U \to U^{\mathbf{p}}(\mathcal{T})$ defined by
\begin{equation}
(u - \Pi^{hp}u, w_{hp}) = 0, \quad \forall w_{hp} \in U^{\mathbf{p}}(\mathcal{T}).
\label{P^hp_dingyi}
\end{equation}
\begin{Theorem}\label{houyanwucha}
	Let $(y,u,z)$ and ($y_{hp},u_{hp},z_{hp}$) be the solution of $\eqref{yijiexitong}$ and $\eqref{lisanxitong}$ respectively. Then $\beta z_{hp}+\lambda u_{hp}\in H^*(\Omega,\mathcal{T})$, we have that
	\begin{equation}
	\|u-u_{hp}\|_{L^2(\Omega)}^2+\|y-y_{hp}\|_{H^1(\Omega)}^2+\|z-z_{hp}\|_{H^1(\Omega)}^2\leq C\eta^2.
	\nonumber
	\end{equation}
\end{Theorem}
\noindent Proof.
To estimate $\|y - y_{hp}\|_{H^1(\Omega)}$ and $\|z - z_{hp}\|_{H^1(\Omega)}$ using the triangle inequality, we first bound the errors $\|y - y(u_{hp})\|_{H^1(\Omega)}$, $\|y(u_{hp}) - y_{hp}\|_{H^1(\Omega)}$, $\|z - z(u_{hp})\|_{H^1(\Omega)}$, and $\|z(u_{hp}) - z_{hp}\|_{H^1(\Omega)}$ individually.\\
We begin by estimating $\|y - y(u_{hp})\|_{H^1(\Omega)}$ and $\|z - z(u_{hp})\|_{H^1(\Omega)}$. By subtracting equations (i)(ii) of \eqref{yijiexitong} from equations (i)(ii) of \eqref{fuzhuxitong}, we obtain:
\begin{align}
a(y(u_{hp}) - y, v) &= (\beta (u_{hp} - u), v)_{\Omega}, \label{first_xiangjian1} \\
a(q, z(u_{hp}) - z) &= (\lambda_\Omega (y(u_{hp}) - y), q)_{\Omega} + (\lambda_\Gamma (y(u_{hp}) - y), q)_{\Gamma}, \label{first_xiangjian2}
\end{align}
Choosing \( v = y(u_{hp}) - y \) in \eqref{first_xiangjian1}, we derive:  
\[
\begin{aligned}
&c\|y(u_{hp}) - y\|_{H^1(\Omega)}^2 \leq a(y(u_{hp}) - y, y(u_{hp}) - y)= (\beta (u_{hp} - u), y(u_{hp}) - y)_{\Omega} \\
&\leq \delta \|u_{hp} - u\|_{L^2(\Omega)} \|y(u_{hp}) - y\|_{L^2(\Omega)}\leq \delta \|u_{hp} - u\|_{L^2(\Omega)} \|y(u_{hp}) - y\|_{H^1(\Omega)}.
\end{aligned}
\]
Canceling the common factor, we obtain:
\begin{equation}
\|y(u_{hp}) - y\|_{H^1(\Omega)} \leq C \|u_{hp} - u\|_{L^2(\Omega)}.
\nonumber
\end{equation}
Analogously, let $ q = z(u_{hp}) - z $ and applying the trace theorem:
\begin{equation}
\begin{aligned}
&c \|z(u_{hp}) - z\|_{H^1(\Omega)}^2 \leq a(z(u_{hp}) - z, z(u_{hp}) - z) \\
&= (\lambda_\Omega (y(u_{hp}) - y), z(u_{hp}) - z)_{\Omega} + (\lambda_\Gamma (y(u_{hp}) - y), z(u_{hp}) - z)_{\Gamma} \\
&\leq \delta \|y(u_{hp}) - y\|_{L^2(\Omega)} \|z(u_{hp}) - z\|_{L^2(\Omega)} + \delta \|y(u_{hp}) - y\|_{L^2(\Gamma)} \|z(u_{hp}) - z\|_{L^2(\Gamma)} \\
&\leq C \|y(u_{hp}) - y\|_{H^1(\Omega)} \|z(u_{hp}) - z\|_{H^1(\Omega)}.
\end{aligned}
\nonumber
\end{equation}
This yields:
\begin{equation}
\|z(u_{hp}) - z\|_{H^1(\Omega)} \leq C \|y(u_{hp}) - y\|_{H^1(\Omega)} \leq C \|u_{hp} - u\|_{L^2(\Omega)}.
\label{4_Theorem_jieguo2}
\end{equation}
Next, we estimate the errors $ \|y(u_{hp}) - y_{hp}\|_{H^1(\Omega)} $ and $\|z(u_{hp}) - z_{hp}\|_{H^1(\Omega)}$. Subtracting (i)(ii) of $\eqref{fuzhuxitong}$ from (i)(ii) of $\eqref{lisanxitong}$, we derive:
\begin{align}
a(y(u_{hp}) - y_{hp}, v_{hp}) &= 0, \quad \forall v_{hp} \in Y^{hp}, \label{second_xiangjian1} \\
a(q_{hp}, z(u_{hp}) - z_{hp}) &= (\lambda_\Omega (y(u_{hp}) - y_{hp}), q_{hp})_{\Omega} + (\lambda_\Gamma (y(u_{hp}) - y_{hp}), q_{hp})_{\Gamma},\nonumber\\
\forall q_{hp} \in Y^{hp}. \label{second_xiangjian2}
\end{align}
Let \( E^y = y(u_{hp}) - y_{hp} \) and \( E_1^y = I_1^{hp} E^y \), where \( I_1^{hp} \) is defined in Lemma \ref{Scott1wenzi}. Then, combining \eqref{second_xiangjian1}, (i) of \eqref{lisanxitong}, (i) of \eqref{fuzhuxitong}, green's formula, and Lemma \ref{Scott1wenzi}, we derive:
\begin{equation}
\begin{aligned}
&c\|y(u_{hp})-y_{hp}\|_{H^1(\Omega)}^2\leq a(E^y,E^y)=a(E^y,E^y-E_1^y)+a(E^y,E_1^y)\\
&=a(E^y,E^y-E_1^y)=a(y(u_{hp})-y_{hp},E^y-E_1^y)\\
&=a(y(u_{hp}),E^y-E_1^y)-\int_{\Omega}\nabla y_{hp}\nabla(E^y-E_1^y)-\int_{\Gamma}\alpha y_{hp}(E^y-E_1^y)\\
&=(\beta u_{hp},E^y-E_1^y)_{\Omega}-\int_{\Omega}\nabla y_{hp}\nabla(E^y-E_1^y)-\int_{\Gamma}\alpha y_{hp}(E^y-E_1^y)\\
&=\int_{\Omega}\beta u_{hp}(E^y-E_1^y)-\int_{\Omega}\nabla y_{hp}\nabla(E^y-E_1^y)-\int_{\Gamma}\alpha y_{hp}(E^y-E_1^y)\\
&=\sum\limits_{\tau\in \mathcal{T}}\left\{\int_{\tau}\Delta y_{hp}(E^y-E_1^y)-\int_{\partial \tau}[\Delta y_{hp}.n_\tau](E^y-E_1^y)\right\}+\int_{\Omega}\beta u_{hp}(E^y-E_1^y)-\int_{\Gamma}\alpha y_{hp}(E^y-E_1^y)\\
&=\sum\limits_{\tau\in \mathcal{T}}\left\{\int_{\tau}(\beta u_{hp}+\Delta y_{hp})(E^y-E_1^y)-\int_{\partial \tau}[\Delta y_{hp}.n_\tau](E^y-E_1^y) \right\}-\sum\limits_{e\in \mathcal{E}_{0}(\mathcal{T})}\int_{e}\alpha y_{hp}(E^y-E_1^y)\\
&=\sum\limits_{\tau\in \mathcal{T}}\int_{\tau}(\beta u_{hp}+\Delta y_{hp})(E^y-E_1^y)-\sum\limits_{e\in \mathcal{E}(\mathcal{T})\backslash\mathcal{E}_{0}(\mathcal{T})}\int_{e}[\nabla y_{hp}.n_e](E^y-E_1^y)\\
&-\sum\limits_{e\in \mathcal{E}_{0}(\mathcal{T})}\int_{e}(\alpha y_{hp}+[\nabla y_{hp}.n_e])(E^y-E_1^y)\\
&\leq C\sum\limits_{\tau\in \mathcal{T}}\|\beta u_{hp}+\Delta y_{hp}\|_{L^2(\tau)}\|E^y-E_1^y\|_{L^2(\tau)}\\
&+C\sum\limits_{e\in \mathcal{E}(\mathcal{T})\backslash\mathcal{E}_{0}(\mathcal{T})}\|[\nabla y_{hp}.n_e]\|_{L^2(e)}\|E^y-E_1^y\|_{L^2(e)}\\
&+C\sum\limits_{e\in \mathcal{E}_{0}(\mathcal{T})}\|\alpha y_{hp}+[\nabla y_{hp}.n_e]\|_{L^2(e)}\|E^y-E_1^y\|_{L^2(e)}\\
&\leq C\sum\limits_{\tau\in \mathcal{T}}\frac{h_\tau^2}{p_\tau^2}\|\beta u_{hp}+\Delta y_{hp}\|_{L^2(\Omega)}^2+C\sum\limits_{e\in \mathcal{E}(\mathcal{T})\backslash\mathcal{E}_{0}(\mathcal{T})}\frac{h_e}{p_e}\|\nabla y_{hp}.n_e\|_{L^2(e)}^2\\
&+C\sum\limits_{e\in \mathcal{E}_{0}(\mathcal{T})}\frac{h_e}{p_e}\|\alpha y_{hp}+[\nabla y_{hp}.n_e]\|_{L^2(e)}^2+\frac{c}{2}\|E^y\|_{H^1(\Omega)}^2,
\end{aligned}
\nonumber
\end{equation}
which meaning, we have the following conclusion,
\begin{equation}
\|y(u_{hp})-y_{hp}\|_{H^1(\Omega)}^2\leq C(\eta_1^2+\eta_2^2+\eta_3^2).
\label{4_Theorem_jieguo3}
\end{equation}
Let $E^z=z(u_{hp})-z_{hp}$ and $E_1^z=I_1^{hp}E^z$. Then it follows from $\eqref{second_xiangjian2}$, item (ii) of $\eqref{lisanxitong}$, green's formula, trace therorm, Lemma \ref{Scott1wenzi},
\begin{align}
&c\|z(u_{hp})-z_{hp}\|_{H^1(\Omega)}^2\leq a(E^z,E^z)=a(z(u_{hp}),E^z)\nonumber-a(z_{hp},E^z-E_1^z)-a(z_{hp},E_1^z)\nonumber\\
&=(\lambda_{\Omega}(y(u_{hp})-y_{\Omega}),E^z)_{\Omega}+(\lambda_{\Gamma}(y(u_{hp})-y_{\Gamma}),E^z)_{\Gamma}-a(z_{hp},E^z-E_1^z)\nonumber\\
&-(\lambda_{\Omega}(y(u_{hp})-y_{\Omega}),E_1^z)_{\Omega}+(\lambda_{\Gamma}(y(u_{hp})-y_{\Gamma}),E_1^z)_{\Gamma}\nonumber\\
&=(\lambda_{\Omega}(y(u_{hp})-y_{hp}),E^z)_{\Omega}+(\lambda_{\Gamma}(y(u_{hp})-y_{hp}),E^z)_{\Gamma}-a(z_{hp},E^z-E_1^z)\nonumber\\
&+(\lambda_{\Omega}(y_{hp}-y_{\Omega}),E^z-E_1^z)_{\Omega}+(\lambda_{\Gamma}(y_{hp}-y_{\Gamma}),E^z-E_1^z)_{\Gamma}\nonumber\\
&=(\lambda_{\Omega}(y(u_{hp})-y_{hp}),E^z)_{\Omega}+(\lambda_{\Gamma}(y(u_{hp})-y_{hp}),E^z)_{\Gamma}\nonumber+(\lambda_{\Omega}(y_{hp}-y_{\Omega}),E^z-E_1^z)_{\Omega}+(\lambda_{\Gamma}(y_{hp}-y_{\Gamma}),E^z-E_1^z)_{\Gamma}\nonumber\\&-\int_{\Gamma}\nabla z_{hp}\nabla(E^z-E_1^z)-\int_{\Gamma}\alpha z_{hp}(E^z-E_1^z)\nonumber\\
&=(\lambda_{\Omega}(y(u_{hp})-y_{hp}),E^z)_{\Omega}+(\lambda_{\Gamma}(y(u_{hp})-y_{hp}),E^z)_{\Gamma}\nonumber+(\lambda_{\Omega}(y_{hp}-y_{\Omega}),E^z-E_1^z)_{\Omega}+(\lambda_{\Gamma}(y_{hp}-y_{\Gamma}),E^z-E_1^z)_{\Gamma}\nonumber\\
&=\sum\limits_{\tau\in \mathcal{T}}\left\{\int_{\tau}\Delta z_{hp}(E^z-E_1^z)-\int_{\partial \tau}[\nabla z_{hp}.n](E^z-E_1^z)\right\}\nonumber-\sum\limits_{e\in \mathcal{E}_{0}(\mathcal{T})}\int_{e}\alpha z_{hp}(E^z-E_1^z)\nonumber\\
&+(\lambda_{\Omega}(y(u_{hp})-y_{hp}),E^z)_{\Omega}+(\lambda_{\Gamma}(y(u_{hp})-y_{hp}),E^z)_{\Gamma}\nonumber+(\lambda_{\Omega}(y_{hp}-y_{\Omega}),E^z-E_1^z)_{\Omega}+(\lambda_{\Gamma}(y_{hp}-y_{\Gamma}),E^z-E_1^z)_{\Gamma}\nonumber\\
&=\sum\limits_{\tau\in \mathcal{T}}\left\{\int_{\tau}(\lambda_{\Omega}(y_{hp}-y_{\Omega})+\Delta z_{hp})(E^z-E_1^z)-\int_{\partial \tau}[\nabla z_{hp}.n](E^z-E_1^z)\right\}\nonumber\\
&+\sum\limits_{e\in \mathcal{E}_{0}(\mathcal{T})}\int_{e}(\lambda_{\Gamma}(y_{hp}-y_{\Gamma})-\alpha z_{hp})(E^z-E_1^z)+(\lambda_{\Omega}(y(u_{hp})-y_{hp}),E^z)_{\Omega}\nonumber+(\lambda_{\Gamma}(y(u_{hp})-y_{hp}),E^z)_{\Gamma}\nonumber\\
&\leq \sum\limits_{\tau\in \mathcal{T}}\|\lambda_{\Omega}(y_{hp}-y_{\Omega})+\Delta z_{hp}\|_{L^2(\tau)}\|E^z-E_1^z\|_{L^2(\tau)}\nonumber+\sum\limits_{e\in \mathcal{E}(\mathcal{T})\backslash\mathcal{E}_{0}(\mathcal{T})}\|[\nabla z_{hp}.n_e]\|_{L^2(e)}\|E^z-E_1^z\|_{L^2(e)}\nonumber\\
&+\sum\limits_{e\in \mathcal{E}_{0}(\mathcal{T})}\|\lambda_{\Gamma}(y_{hp}-y_{\Gamma})-\alpha z_{hp}-[\nabla z_{hp}.n_e]\|_{L^2(e)}\|E^z-E_1^z\|_{L^2(e)}\nonumber\\
&+\|\lambda_{\Omega}(y(u_{hp})-y_{hp})\|_{L^2(\Omega)}\|E^z\|_{L^2(\Omega)}+\|\lambda_{\Gamma}(y(u_{hp})-y_{hp})\|_{L^2(\Omega)}\|E^z\|_{L^2(\Gamma)}\nonumber\\
&\leq C\sum\limits_{\tau\in \mathcal{T}}\frac{h_\tau^2}{p_\tau^2}\|\lambda_{\Omega}(y_{hp}-y_{\Omega})+\Delta z_{hp}\|_{L^2(\tau)}^2+C\sum\limits_{e\in \mathcal{E}(\mathcal{T})\backslash\mathcal{E}_{0}(\mathcal{T})}\frac{h_e}{p_e}\|[\nabla z_{hp}.n_e]\|_{L^2(e)}^2\nonumber\\
&+C\sum\limits_{e\in \mathcal{E}_{0}(\mathcal{T})}\frac{h_e}{p_e}\|\lambda_{\Gamma}(y_{hp}-y_{\Gamma})-\alpha z_{hp}-[\nabla p_{hp}.n_e]\|_{L^2(e)}^2\nonumber+C\|y(u_{hp})-y_{hp}\|_{H^1(\Omega)}^2+\frac{c}{2}\|E^z\|_{H^1(\Omega)}^2.
\nonumber
\end{align}
\begin{equation}
\|z(u_{hp})-z_{hp}\|_{H^1(\Omega)}^2\leq C(\eta_4^2+\eta_5^2+\eta_6^2)+C\|y(u_{hp})-y_{hp}\|_{H^1(\Omega)}^2.
\label{4_Theorem_jieguo4}
\end{equation}
We estimate the error in the final part $\|u-u_{hp}\|_{L^2(\Omega)}^2$,\\
\begin{equation}
c\|u-u_{hp}\|_{L^2(\Omega)}^2\leq (\lambda u-\lambda u_{hp},u-u_{hp})=(\lambda u,u-u_{hp})-(\lambda u_{hp},u-u_{hp}).\\
\nonumber
\end{equation}
Using $\eqref{lambda_u_xiaoyu1}$ and $\eqref{lambda_u_xiaoyu2}$, $\forall v_{hp}\in U^{hp}$, we have
\begin{equation}
\begin{aligned}
&(\lambda u,u-u_{hp})-(\lambda u_{hp},u-u_{hp})\leq(-\beta z,u-u_{hp})+(\beta z_{hp},u-u_{hp})+(\beta z_{hp}+\lambda u_{hp},v_{hp}-u)\\
&=(\beta(z_{hp}-z),u-u_{hp})+(\beta z_{hp}+\lambda u_{hp},v_{hp}-u)\\
&\leq(\beta z_{hp}+\lambda u_{hp},v_{hp}-u)+(\beta(z_{hp}-z(u_{hp})),u-u_{hp})+(\beta(z(u_{hp})-z),u-u_{hp})\\
&\leq(\beta z_{hp}+\lambda u_{hp},v_{hp}-u)+(\beta(z_{hp}-z(u_{hp})),u-u_{hp})-((z(u_{hp})-z),\beta(u_{hp}-u)).
\end{aligned}
\label{lamdba_u-u_hp}
\end{equation}
We have
\begin{equation}
\begin{aligned}
&(\beta(u_{hp}-u),z(u_{hp})-z)\\
&=(\lambda_{\Omega}(y(u_{hp})-y),y(u_{hp})-y)_{\Omega}+(\lambda_{\Gamma}(y(u_{hp})-y),y(u_{hp})-y)_{\Gamma}\geq C(\|y(u_{hp})-y\|_{L^2(\Omega)}^2+\|y(u_{hp})-y\|_{L^2(\Gamma)}^2)\geq 0.
\nonumber
\end{aligned}
\end{equation}
This implies  
\[
-(z(u_{hp}) - z, \beta(u_{hp} - u)) \leq 0.
\]  
Rewriting \eqref{lamdba_u-u_hp}, we obtain
\begin{equation}
c\|u-u_{hp}\|_{L^2(\Omega)}^2\leq(\beta z_{hp}+\lambda u_{hp},v_{hp}-u)+(\beta(z_{hp}-z(u_{hp})),u-u_{hp})
\label{lamdba_u-u_hp_huajian1}
\end{equation}
Next, we analyze the two terms on the right-hand side of \eqref{lamdba_u-u_hp_huajian1} separately. Following the techniques in \cite{chen2011posteriori}, we analogously define
\begin{equation}
G\lvert_{\tau}=\beta z_{hp}+\lambda u_{hp}\quad \widehat{G}\lvert_{\hat{\tau}}=\widehat{G\lvert_{\tau}}=(\beta z_{hp}+\lambda u_{hp})\circ F_{\tau}\in H^1(\hat{\tau})\quad \overline{\widehat{G}\lvert_{\hat{\tau}}}=\int_{\hat{\tau}}\widehat{G}\lvert_{\hat{\tau}}/\int_{\hat{\tau}}1.
\nonumber
\end{equation}
Let \( v_{hp} = \mathscr{R}^{hp}u \in U^{hp} \subset U^{\mathbf{p}}(\mathcal{T}) \). From \eqref{fuzhuxitong}, the definition of \( I_U^{hp} \), and \eqref{P^hp_dingyi},
\begin{equation}
\begin{aligned}
&(\beta z_{hp}+\lambda u_{hp},v_{hp}-u)=(\beta z_{hp}+\lambda u_{hp},\Pi^{hp}u-u)\\
&=(\beta z_{hp}+\lambda u_{hp}-I_U^{hp}(\beta z_{hp}+\lambda u),\Pi^{hp}u-u)+(I_U^{hp}(\beta z_{hp}+\lambda u),\Pi^{hp}u-u)\\
&=(\beta z_{hp}+\lambda u_{hp}-I_U^{hp}(\beta z_{hp}+\lambda u),\Pi^{hp}u-u).
\end{aligned}
\nonumber
\end{equation}
It follows from $\eqref{P^hp_dingyi}$,
\begin{equation}
\begin{aligned}
&\|u-\Pi^{hp}u\|\|u-v_{hp}\|\geq (u-\Pi^{hp}u,u-v_{hp})=(u-\Pi^{hp}u,u-\Pi^{hp}u)+(u-\Pi^{hp}u,\Pi^{hp}u-v_{hp})\geq\|u-\Pi^{hp}u\|^2.
\nonumber
\end{aligned}
\end{equation}
Then we have $\|u-v_{hp}\|\geq \|u-\Pi^{hp}u\|$ , which implies 
\begin{equation}
\|u-u_{hp}\|\geq \|u-\Pi^{hp}u\|.
\label{u_First_of_right_2}
\end{equation}
It follows from Lemma \ref{Ihp1wenzi}, $\eqref{P^hp_dingyi}$, Cauchy–Schwarz inequality, Poincaré inequality, and $\eqref{u_First_of_right_2}$,
\begin{equation}
\begin{aligned}
&(\beta z_{hp}+\lambda u_{hp},v_{hp}-u)\\
&=\sum\limits_{\tau\in \mathcal{T}}(\beta z_{hp}+\lambda u_{hp}-\pi_{p_{\tau}}^{h_{\tau}}(\beta z_{hp}+\lambda u_{hp}),\Pi^{hp}u-u)_{L^2(\tau)}\\
&=\sum\limits_{\tau\in \mathcal{T}}(\beta z_{hp}+\lambda u_{hp}-\overline{\widehat{G}\lvert_{\hat{\tau}}}-\pi_{p_{\tau}}^{h_{\tau}}(\beta z_{hp}+\lambda u_{hp}-\overline{\widehat{G}\lvert_{\hat{\tau}}}),\Pi^{hp}u-u)_{L^2(\tau)}\\
&+\sum\limits_{\tau\in \mathcal{T}}(\pi_{p_{\tau}}^{h_{\tau}}(\beta z_{hp}+\lambda u_{hp}-\overline{\widehat{G}\lvert_{\hat{\tau}}})-\pi_{p_{\tau}}^{h_{\tau}}(\lambda u_{hp}+\beta z_{hp})+\overline{\widehat{G}\lvert_{\hat{\tau}}},\Pi^{hp}u-u)_{L^2(\tau)}\\
&=\sum\limits_{\tau\in \mathcal{T}}(\beta z_{hp}+\lambda u_{hp}-\overline{\widehat{G}\lvert_{\hat{\tau}}}-\pi_{p_{\tau}}^{h_{\tau}}(\beta z_{hp}+\lambda u_{hp}-\overline{\widehat{G}\lvert_{\hat{\tau}}}),\Pi^{hp}u-u)_{L^2(\tau)}\\
&\leq C\sum\limits_{\tau\in \mathcal{T}}\frac{h_\tau}{p_\tau}\|\lambda u_{hp}+\beta z_{hp}-\overline{\widehat{G}\lvert_{\hat{\tau}}}\|_{H^1(\tau)}\|\Pi^{hp}u-u\|_{L^2(\tau)}\\
&\leq C\sum\limits_{\tau\in \mathcal{T}}\frac{h_{\tau}^2}{p_{\tau}^2}\|\lambda u_{hp}+\beta z_{hp}-\overline{\widehat{G}\lvert_{\hat{\tau}}}\|_{H^1(\tau)}^2+\frac{c}{4}\|\Pi^{hp}u-u\|_{L^2(\Omega)}^2\\
&\leq C\sum\limits_{\tau\in \mathcal{T}}\frac{h_{\tau}^2}{p_{\tau}^2}\|\nabla(\lambda u_{hp}+\beta z_{hp})\|_{L^2(\tau)}^2+\frac{c}{4}\|\Pi^{hp}u-u\|_{L^2(\Omega)}^2\\
&\leq C\sum\limits_{\tau\in \mathcal{T}}\frac{h_{\tau}^2}{p_{\tau}^2}\|\nabla(\lambda u_{hp}+\beta z_{hp})\|_{L^2(\tau)}^2+\frac{c}{4}\|u-u_{hp}\|_{L^2(\Omega)}^2.
\end{aligned}
\label{u_First_of_right_3}
\end{equation}
Next, we estimate the second part of $\eqref{lamdba_u-u_hp_huajian1}$. It follows from Cauchy-Schwarz inequality and Young inequality,
\begin{equation}
\begin{aligned}
&(\beta(z_{hp}-z(u_{hp})),u-u_{hp})\leq \|\beta(z_{hp}-z(u_{hp})\|_{L^2(\Omega)}\|u-u_{hp}\|_{L^2(\Omega)}\\
&\leq C\|z_{hp}-z(u_{hp})\|_{H^1(\Omega)}^2+\frac{c}{2}\|u-u_{hp}\|_{L^2(\Omega)}^2.
\end{aligned}
\label{u_Second_of_right_1}
\end{equation}
Combining \eqref{lamdba_u-u_hp_huajian1}, \eqref{u_First_of_right_3}, and \eqref{u_Second_of_right_1}, we obtain
\begin{equation}
\begin{aligned}
&c\|u-u_{hp}\|_{L^2(\Omega)}^2\leq(\beta z_{hp}+\lambda u_{hp},v_{hp}-u)+(\beta(z_{hp}-z(u_{hp})),u-u_{hp})\\
&\leq C\sum\limits_{\tau\in \mathcal{T}}\frac{h_{\tau}^2}{p_{\tau}^2}\|\nabla(\lambda u_{hp}+\beta z_{hp})\|_{L^2(\tau)}^2+C\|z_{hp}-z(u_{hp})\|_{H^1(\Omega)}^2+\frac{3c}{4}\|u-u_{hp}\|_{L^2(\Omega)}^2.\\
\nonumber
\end{aligned} 
\end{equation}
Moving the term \(\|u - u_{hp}\|_{L^2(\Omega)}^2\) to the right-hand side, we obtain
\begin{equation}
\begin{aligned}
\|u-u_{hp}\|_{L^2(\Omega)}^2&\leq C\sum\limits_{\tau\in \mathcal{T}}\frac{h_{\tau}^2}{p_{\tau}^2}\|\nabla(\lambda u_{hp}+\beta z_{hp})\|_{L^2(\tau)}^2+C\|z_{hp}-z(u_{hp})\|_{H^1(\Omega)}^2\\
&=C\eta_7^2+C\|z_{hp}-z(u_{hp})\|_{H^1(\Omega)}^2.
\end{aligned}
\label{4_Theorem_jieguo5}
\end{equation}
Combining \eqref{4_Theorem_jieguo2}, \eqref{4_Theorem_jieguo3}, \eqref{4_Theorem_jieguo4}, and \eqref{4_Theorem_jieguo5},
\begin{equation}
\begin{aligned}
&\|u-u_{hp}\|_{L^2(\Omega)}^2+\|y-y_{hp}\|_{H^1(\Omega)}^2+\|z-z_{hp}\|_{H^1(\Omega)}^2\\
&\leq\|u-u_{hp}\|_{L^2(\Omega)}^2+2(\|y-y(u_{hp})\|_{H^1(\Omega)}^2+\|z-z(u_{hp})\|_{H^1(\Omega)}^2)\\
&+2(\|y(u_{hp})-y_{hp}\|_{H^1(\Omega)}^2+\|z(u_{hp})-z_{hp}\|_{H^1(\Omega)}^2)\\
&\leq C\|u-u_{hp}\|_{L^2(\Omega)}^2+2(\|y(u_{hp})-y_{hp}\|_{H^1(\Omega)}^2+\|z(u_{hp})-z_{hp}\|_{H^1(\Omega)}^2)\\
&\leq C\eta^2.
\end{aligned}
\nonumber
\end{equation}
Then proof is complete. $\hfill\square$
\subsection{Numerical algorithm}
An iterative method that has been proven effective \cite{MR3559582} is provided in this section, which derived from the first-order optimality conditions, follows the logical update sequence:
\begin{equation}
u_{hp}^n \rightarrow y_{hp}^n \rightarrow z_{hp}^n \rightarrow u_{hp}^{n+1} \rightarrow \cdots\nonumber
\end{equation}
Step 1: solve for the adjoint variable $z$\\
Given the current iterates $u_{hp}^n$ and $y_{hp}^{n}$, the adjoint variable $z_{hp}^{n}$ is obtained by solving the following variational formulation:
\begin{equation}
a(q_{hp}, z_{hp}^{n}) = (\lambda_{\Omega}(y_{hp}^{n} - y_{\Omega}), q_{hp})_{\Omega} + (\lambda_{\Gamma}(y_{hp}^{n} - y_{\Gamma}), q_{hp})_{\Gamma}, \quad \forall q_{hp} \in Y^{hp}.\nonumber
\end{equation}
The solution $z_{hp}^{n}$ is then used to update the control variable in the subsequent step.\\
Step 2: control update (projected gradient descent)\\
This unconstrained solution is then projected onto the feasible set $U_{a}=\{u>u_a\}$ via the projection method, resulting in the updated control:
\begin{equation}
u_{hp}^{n+1}=\max\{u_a,-\frac{\beta}{\lambda}z^{n}\}\nonumber
\end{equation}
Given the updated control $u_{hp}^n$, the new state $y_{hp}^n$ is computed by solving the variational formulation:
\begin{equation}
a(y_{hp}^{n+1},v_{hp})=(\beta u_{hp}^{n+1},v_{hp})_\Omega\quad\forall v_{hp}\in Y^{hp}\nonumber
\end{equation}
This problem is identical in form to the initial state equation but is driven by the new control $u_{hp}^{n+1}$. The solution $y_{hp}^{n+1}$ then provides the updated state for the subsequent adjoint equation in the next iteration.\\

%
\begin{algorithm}[H]
	\caption{Iteration method}
	\label{algorithm}
	\begin{algorithmic}[1]
		\Require Given constant $tol, u_a\in\mathbb{R}$, select an initial values $u_{hp}^{0}\in U^{hp}$, seek $y_{hp}^{0}\in Y^{hp}$ such that
		\begin{equation}
			a(y_{hp}^{0},v_{hp})=(\beta u_{hp}^{0},v_{hp})_\Omega,\quad\forall v_{hp}\in Y^{hp}.\nonumber
		\end{equation}
		\State Seek $z_{hp}^{n}\in Y^{hp}$ such that
		\begin{equation}
		a(q_{hp},z_{hp}^{n})=(\lambda_{\Omega}(y_{hp}^{n}-y_{\Omega}),q_{hp})_{\Omega}+(\lambda_{\Gamma}(y_{hp}^{n}-y_{\Gamma}),q_{hp})_{\Gamma},  \quad \forall q_{hp}\in Y^{hp}.\nonumber
		\end{equation}
		\State Seek $y_{hp}^{n+1}\in Y^{hp}$ such that
        \begin{equation}
        	a(y_{hp}^{n+1},v_{hp})=(\beta u_{hp}^{n+1},v_{hp})_\Omega,\quad\forall v_{hp}\in Y^{hp}.\nonumber
        \end{equation}
		\State {Stop if stopping criterion $\|z_{hp}^{n+1}-z_{hp}^{n}\|_{H^1(\Omega)}\leq tol$ is satisfied. Otherwise set $n=n+1$ and then back to Step 1.}
	\end{algorithmic}
\end{algorithm}

\section{Numerical experiments}\label{sec5}
We validate the theoretical results by a numerical experiment. The error analysis of our triangular element solutions will be performed against a high-precision reference solution computed using the standard finite element method on a fine mesh of 2500 rectangular elements. \\
\textbf{Example 1} We consider problem \eqref{eq1}-\eqref{eq4} with $\Omega=(0,1)^2\in \mathbb{R}^2$ and $u_a=0$. The target functions and all parameters are given by
\begin{flalign*}
&\lambda = 1/2,\quad\lambda_{\Omega} = \lambda_{\Gamma} = \beta = \alpha = 1,\nonumber\\
&y_{\Omega} = 10x_1 x_2\sin(\pi x_1)\sin(\pi x_2),\quad y_{\Gamma} = y_{\Omega}\lvert_{\Gamma}.\nonumber
\end{flalign*}
In Table \ref{table1}, with the mesh fixed at $N = 64$ elements, the quantities $\|u-u_{hp}\|_{L^2(\Omega)},\|y-y_{hp}\|_{L^2(\Omega)},\|y-y_{hp}\|_{H^1(\Omega)},\|z-z_{hp}\|_{L^2(\Omega)},\|z-z_{hp}\|_{H^1(\Omega)}$ are listed for polynomial degrees 1, 2 Table \ref{table2}, keeps $\mathbf{p} = 2$ fixed and displays the same error norms as the number of elements varies over $N = 16, 64$.\\
\begin{figure}[htbp]
	\centering
	\includegraphics[width=0.8\textwidth]{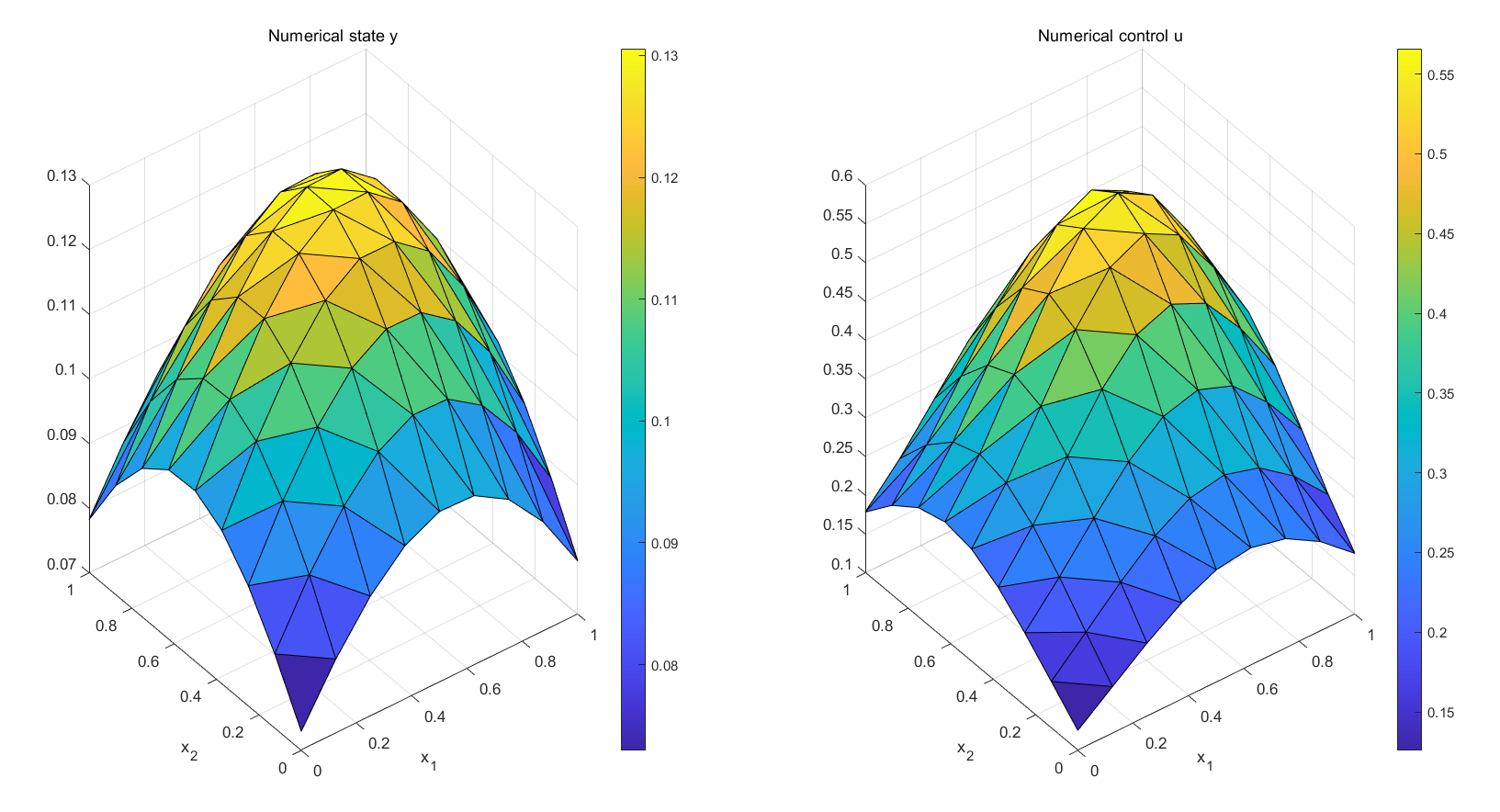}
	\caption{Numerical solution for $N=64$ and $\mathbf{p}=2$}
	\label{fig:example2}
\end{figure}
\begin{figure}[t]
	\centering
	\includegraphics[width=1\textwidth]{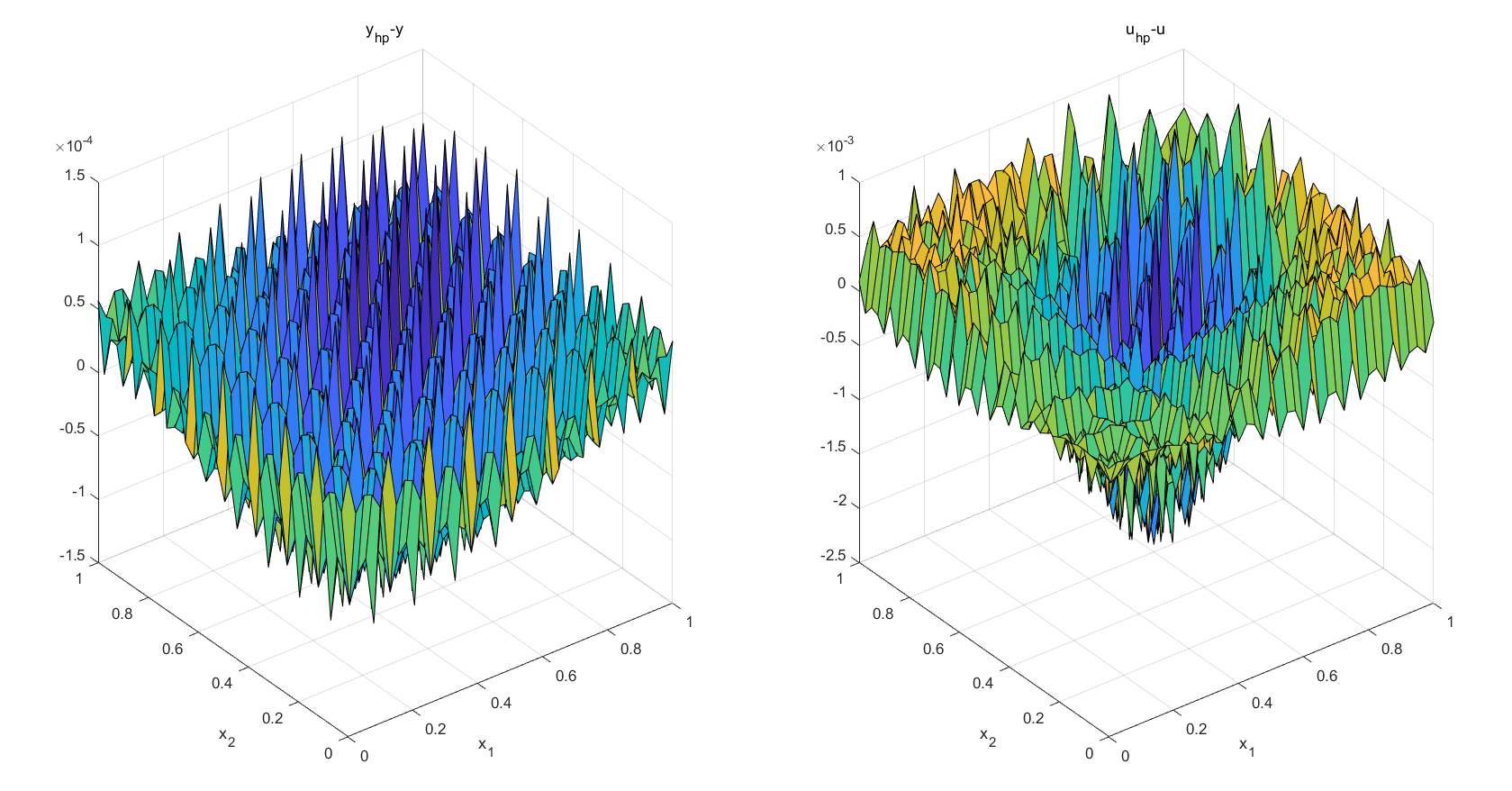}
	\caption{Error for $N=64$ and $\mathbf{p}=2$}
	\label{fig:example2_error}
\end{figure}
\begin{table}[h]
	\centering 
	\caption{Errors vs. polynomial order at number of elements = 64} 
	\label{table1} 
	\setlength{\tabcolsep}{2pt} 
	\renewcommand{\arraystretch}{0.6} 
	\begin{tabular}{l r r r r r}
		\toprule 
		
		$\boldsymbol{p}$ &$\|u-u_{hp}\|_{L^2(\Omega)}$ & $\|y-y_{hp}\|_{L^2(\Omega)}$ & $\|y-y_{hp}\|_{H^1(\Omega)}$ & $\|z-z_{hp}\|_{L^2(\Omega)}$ & $\|z-z_{hp}\|_{H^1(\Omega)}$ \\
		\midrule 
		
		1 & 2.567562e-02 & 1.708764e-02 & 2.867192e-02 & 1.283781e-02 & 5.966962e-02 \\
		
		2 & 6.398634e-04 & 3.939740e-05 & 3.575059e-03 & 3.199326e-04 & 1.476286e-02 \\
		
		\bottomrule 
	\end{tabular}
    \renewcommand{\arraystretch}{1} 
\end{table}
\begin{table}[h]
	\centering 
	\caption{Errors vs. number of elements at polynomial order = 2} 
	\label{table2} 
	\setlength{\tabcolsep}{2pt} 
	\renewcommand{\arraystretch}{0.6} 
	\begin{tabular}{l r r r r r}
		\toprule 
		
		$N$ &$\|u-u_{hp}\|_{L^2(\Omega)}$ & $\|y-y_{hp}\|_{L^2(\Omega)}$ & $\|y-y_{hp}\|_{H^1(\Omega)}$ & $\|z-z_{hp}\|_{L^2(\Omega)}$ & $\|z-z_{hp}\|_{H^1(\Omega)}$ \\
		\midrule 
		
		16 & 3.767134e-03 & 4.946770e-04 & 1.014658e-02 & 1.883568e-03 & 4.141525e-02 \\
		
		64 & 6.398634e-04 & 3.939740e-05 & 3.575059e-03 & 3.199326e-04 & 1.476286e-02 \\
		
		\bottomrule 
	\end{tabular}
    \renewcommand{\arraystretch}{1} 
\end{table}
\begin{table}[h]
	\centering 
	\small
	\caption{Right-hand side of the posteriori error estimator} 
	\label{table3} 
	\resizebox{\textwidth}{!}{ 
		\begin{tabular}{c c c c c c c c} 
			\toprule 
			
			$\eta_1^2$ & $\eta_2^2$ & $\eta_3^2$ & $\eta_4^2$ & $\eta_5^2$ & $\eta_6^2$ & $\eta_7^2$ & $\eta^2$ \\
			\midrule 
			
			9.928605e-04 & 4.456117e-04 & 4.507555e-03 & 1.180462e-02 & 6.448955e-03 & 3.532674e-02 & 1.597850e-23 & 5.952634e-02 \\
			\bottomrule 
		\end{tabular}
	}
\end{table}

\begin{figure}[b]
	\centering
	\includegraphics[width=0.8\textwidth]{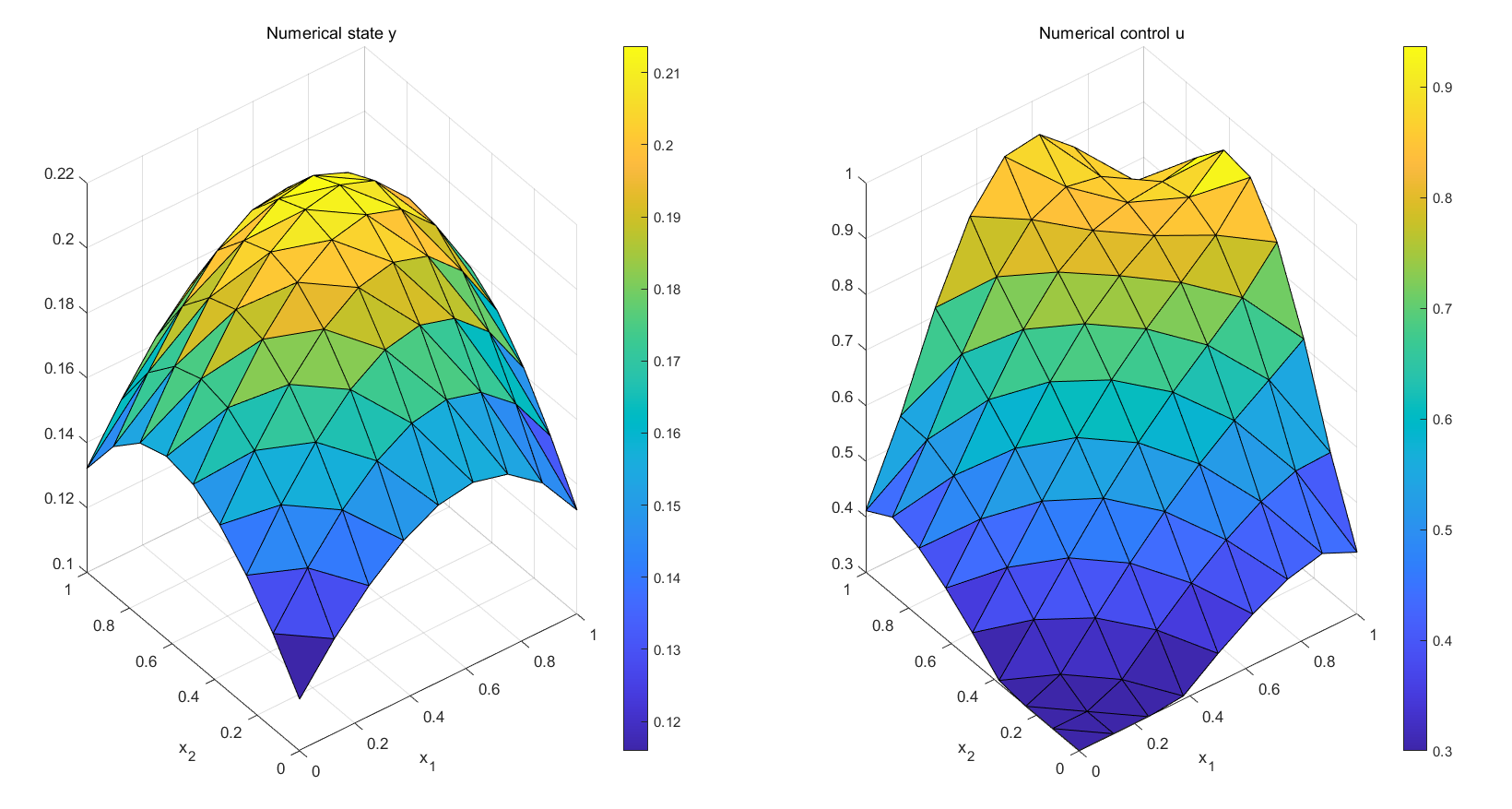}
	\caption{Numerical solution for $N=64$ and $\mathbf{p}=2$}
	\label{fig:example3}
\end{figure}

\begin{figure}[t]
	\centering
	\includegraphics[width=1\textwidth]{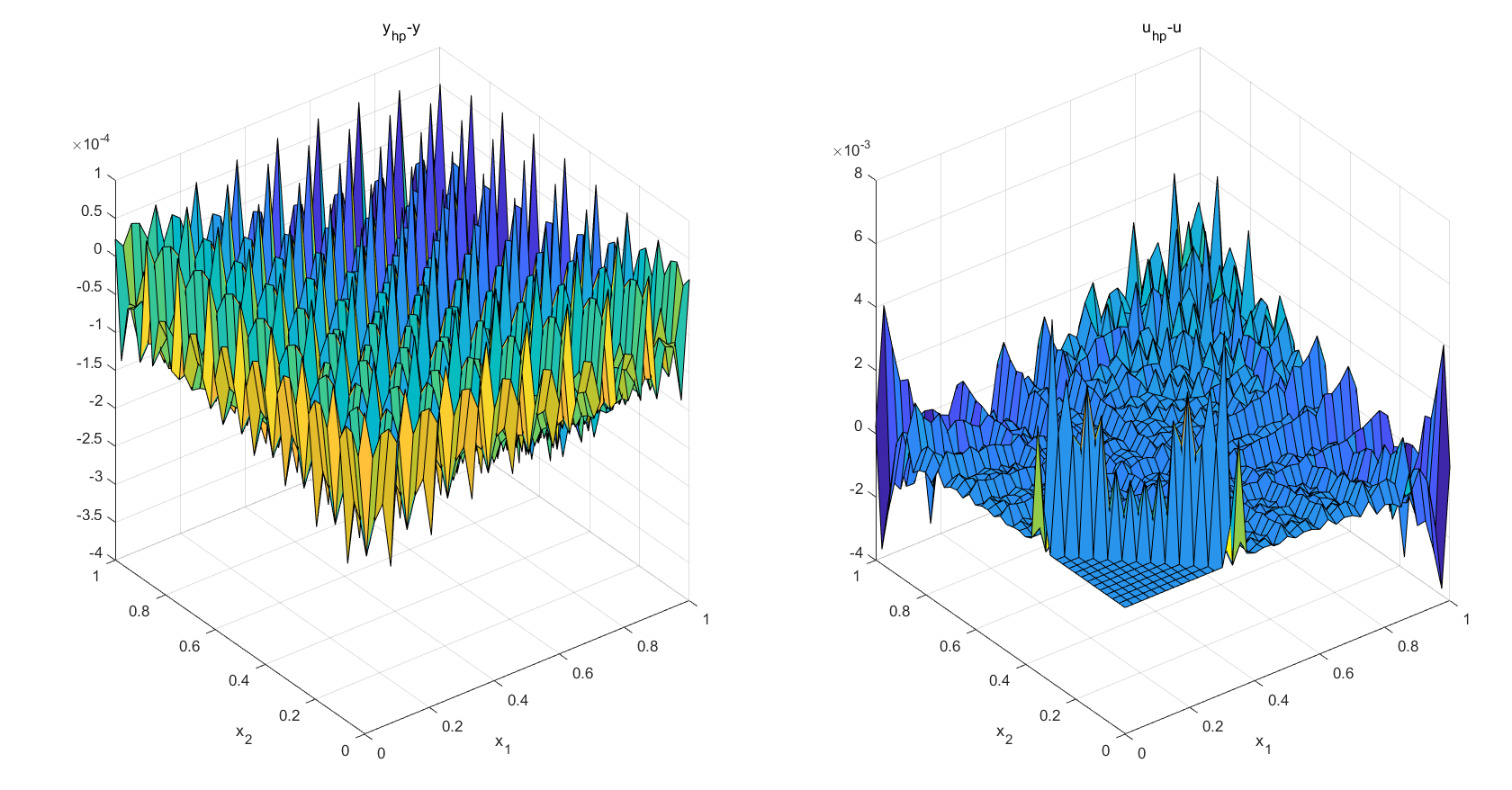}
	\caption{Error for $N=64$ and $\mathbf{p}=2$}
	\label{fig:example3_error}
\end{figure}

\noindent Under the high-resolution setting of $\mathbf{p}=2$ and a fine mesh, Figure \ref{fig:example2} displays the shapes of the numerical solutions $y$ and $u$ across their domain. Figure \ref{fig:example2_error} displays the shapes of errors.\\
Both Table \ref{table1} and Table \ref{table2} demonstrate the convergence of the method. Table \ref{table1} indicates that increasing $\mathbf{p}$ from 1 to 2 yields a reduction in the $L^2$ errors by 2-3 orders of magnitude for $u$, $y$, and $z$, and a reduction in the $H^1$ error by about one order for $y$. Table \ref{table2} shows that refining the mesh from $N=16$ to $N=64$ reduces the $L^2$ errors of $u$, $y$, and $z$ by an order of magnitude, with a corresponding decrease in their $H^1$ errors.\\
Setting $\mathbf{p}=2$, $N=64$, the computed posteriori error estimators are in Talbe \ref{table3}
\noindent The values listed in the second row of Table \ref{table2}, when compared to $\eta^2$, validate the result of Theorem \ref{houyanwucha}. \\
\textbf{Example 2} Set $\Omega=(0,1)^2\in \mathbb{R}^2$ and $u_a=0.3$. The target functions and all parameters are given by
\begin{flalign*}
&\lambda = 1/2,\quad\lambda_{\Omega} = \lambda_{\Gamma} = \beta = \alpha = 1,\nonumber\\
&y_{\Omega} = x_1\sin(\pi x_2)+x_2\sin(\pi x_1),\quad y_{\Gamma} = y_{\Omega}\lvert_{\Gamma}.\nonumber
\end{flalign*}
\begin{table}[h]
	\centering 
	\caption{Errors vs. polynomial order at number of elements = 64} 
	\label{table4} 
	\setlength{\tabcolsep}{2pt} 
	\renewcommand{\arraystretch}{0.6} 
	\begin{tabular}{l r r r r r}
		\toprule 
		
		$\boldsymbol{p}$ &$\|u-u_{hp}\|_{L^2(\Omega)}$ & $\|y-y_{hp}\|_{L^2(\Omega)}$ & $\|y-y_{hp}\|_{H^1(\Omega)}$ & $\|z-z_{hp}\|_{L^2(\Omega)}$ & $\|z-z_{hp}\|_{H^1(\Omega)}$ \\
		\midrule 
		
		1 & 1.385394e-02 & 2.141509e-02 & 4.417284e-02 & 6.906800e-03 & 6.524066e-02 \\
		
		2 & 6.878369e-04 & 1.384393e-04 & 5.875335e-03 & 2.755939e-04 & 1.891997e-02 \\
		
		\bottomrule 
	\end{tabular}
    \renewcommand{\arraystretch}{1} 
\end{table}
\begin{table}[h]
	\centering 
	\caption{Errors vs. number of elements at polynomial order = 2} 
	\label{table5} 
	\setlength{\tabcolsep}{2pt} 
	\renewcommand{\arraystretch}{0.6} 
	\begin{tabular}{l r r r r r}
		\toprule 
		
		$N$ &$\|u-u_{hp}\|_{L^2(\Omega)}$ & $\|y-y_{hp}\|_{L^2(\Omega)}$ & $\|y-y_{hp}\|_{H^1(\Omega)}$ & $\|z-z_{hp}\|_{L^2(\Omega)}$ & $\|z-z_{hp}\|_{H^1(\Omega)}$ \\
		\midrule 
		
		16 & 2.991602e-03 & 7.801649e-04 & 1.665434e-02 & 1.304290e-03 & 5.413024e-02 \\
		
		64 & 6.878369e-04 & 1.384393e-04 & 5.875335e-03 & 2.755939e-04 & 1.891997e-02 \\
		
		\bottomrule 
	\end{tabular}
    \renewcommand{\arraystretch}{1} 
\end{table}
\begin{table}[h]
	\centering 
	\caption{Right-hand side of the posteriori error estimator} 
	\label{table6} 
	\resizebox{\textwidth}{!}{ 
		\begin{tabular}{c c c c c c c c} 
			\toprule 
			
			$\eta_1^2$ & $\eta_2^2$ & $\eta_3^2$ & $\eta_4^2$ & $\eta_5^2$ & $\eta_6^2$ & $\eta_7^2$ & $\eta^2$ \\
			\midrule 
			
			2.797200e-03 & 1.232243e-03 & 1.270797e-02 & 1.218249e-02 & 7.707540e-03 & 3.566816e-02 & 4.502413e-05 & 7.234063e-02 \\
			\bottomrule 
		\end{tabular}
	}
\end{table}
\noindent In this example, we consider a problem with an active control constraint. Figure \ref{fig:example3} and Figure \ref{fig:example3_error} display the numerical solution and the error distribution, respectively. Due to the lower bound constraint $U_a$, the control variable $u$ vanishes in the vicinity of (0,0).\\
Tables \ref{table4} and \ref{table5} demonstrate the convergence of the method, while Tables \ref{table5} and \ref{table6} confirm the effectiveness of the posteriori error estimators.

\clearpage 

\section{Conclusion}\label{sec6}
In this paper, we investigated a priori and a posteriori error analysis of the $hp$-FEM for elliptic optimal control problems with Robin boundary condition. Clément-type interpolation operators and auxiliary systems are employed to derive error bounds for the control, state, and adjoint state variables in the $H^1$ and $L^2$-norms under continuous optimality conditions. For the a posteriori analysis, residual-based error estimators were derived using Scott–Zhang-type quasi interpolation and coupled state-control approximations, providing computable upper bounds for the discretization errors.The theoretical framework presented here extends the applicability of $hp$-FEM to Robin boundary control problems. Numerical results support the theoretical findings.

\bibliographystyle{elsarticle-num}
\bibliography{REF2}

\end{document}